\documentclass{article}

\usepackage{amssymb,amscd,amsthm,latexsym}
\usepackage[all]{xy}

\newtheorem{defi}{Definition}[section]
\newtheorem{theo}{Theorem}[section]
\newtheorem{lemma}{Lemma}[section]
\newtheorem{prop}{Proposition}[section]
\newtheorem{coro}{Corollary}[section]

\newtheorem{ex}{Example}[section]
\def\into{ \rightarrowtail }
\def\onto{ \twoheadrightarrow }

\newdir{ >}{{}*!/-8pt/@{>}}
\def\EE{ \mathbb{E} }
\def\FF{ \mathbb{F} }
\def\CC{ \mathbb{C} }

\begin{document}

\author{Dominique Bourn}

\title{Normalizers in the non-pointed context:\\
a weak case of extremal decomposition}

\date{}

\maketitle

\begin{abstract}
The aim of this work is to point out a strong structural phenomenon hidden behind the existence of normalizers through the investigation of this property in the non-pointed context: given any category $\mathbb E$, a certain property of the fibration of points $\P_{\mathbb{E}} \colon Pt(\mathbb E) \to \mathbb{E}$ guarentees the existence of normalizers. This property becomes a characterization of this existence when $\mathbb E$ is quasi-pointed and protomodular. This property is also showed to be equivalent to a property of the category $Grd\mathbb E$ of internal groupoids in $\mathbb E$ which is a kind of opposite, for the monomorphic internal functors, of the comprehensive factorization.
\end{abstract}

\bigskip

\noindent\textbf{subjclass:} 18A05, 18B99, 18E13, 08C05, 08A30, 08A99.\\
\textbf{keywords:} normal subobject, normalizer, Mal'tsev and protomodular categories, internal categories and groupoids, comprehensive factorization, non-pointed additive categories.

\section*{Introduction}

The first place where the question of the existence of normalizers was investigated in a conceptual way, namely outside specific contexts as groups, rings or Lie algebras, but more generally inside any semi-abelian category (which is a pointed context) is \cite{Gr}. Modulo a slight shifting in the requirement of the involved universal property, it was showed in \cite{BG2} that, in the pointed protomodular context, the existence of normalizers is unexpectedly equivalent to a much larger phenomenon involving the split exact sequences and that it has two heavy structural consequences, namely that the ground $\mathbb C$ is \emph{action accessible} in the sense of \cite {BJ} and \emph{fiberwise algebraically cartesian closed} in the sense of \cite{BG}. 

The notion of normal subobject having a plain meaning in a non-pointed context, the notion of normalizer is straightforward. The aim of this work was first to investigate whether there was, in a non-pointed context, a condition which characterizes the existence of normalizers, or in other words to transfer the pointed characteristic condition of \cite{BG2} to a non-pointed one. Actually we do better, introducing two  equivalent conditions which, here again, are far from being expected (even in the pointed context and even in the categories $Gp$ of groups and $K$-$Lie$ of Lie algebras on a field $K$) and pointing out much larger phenomenons:\\
1) the existence of a universal decomposition for the monomorphisms between split epimorphisms in $\mathbb E$; namely, any monomorphism $(y,x):(f's')\rightarrowtail (f,s)$  produces a universal dotted decomposition, as in the following left hand side diagram, where the left hand side part is a pullback:
$$
\xymatrix@=20pt{
	{X'\;} \ar[d]_{f'} \ar@<2ex>@{>->}[rr]^{x} \ar@{>.>}[r]_{\bar u} & {\bar X\;} \ar[d]_{\bar f} \ar@{>.>}[r]_{\bar w}  & X \ar[d]_{f} && {U_1\;}\ar@{>.>}[r]_{u_1} \ar@{>->}@<2ex>[rr]^{v_1} \ar@<1ex>[d]^{d_1}\ar@<-1ex>[d]_{d_0} & {X_1\;} \ar@<1ex>[d]^{d_1}\ar@<-1ex>[d]_{d_0} \ar@{>.>}[r]_{w_1} &  T_1 \ar@<1ex>[d]^{d_1}\ar@<-1ex>[d]_{d_0}\\
	{Y\;} \ar@<-1ex>[u]_{s'} \ar@<-2ex>@{>->}[rr]_{y} \ar@{>.>}[r]^{u} & {\bar Y\;} \ar@{>.>}[r]^{w} \ar@<-1ex>[u]_{\bar s} & Y \ar@<-1ex>[u]_{s} && {U_0\;} \ar@{>.>}[r]^{u_0} \ar@{>->}@<-2ex>[rr]_{v_0} \ar[u] & {X_0\;} \ar[u] \ar@{>.>}[r]^{w_0}  & T_0 \ar[u]
}
$$
2) the existence of a universal decomposition in the category $Grd\mathbb E$ of internal groupoids in $\mathbb E$ which is a kind of opposite, for the monomorphic functors, of the comprehensive factorization of \cite{SW}, \cite{B-2}; namely, any monomorphic functor  $(v_0,v_1):\underline U_1\rightarrowtail \underline T_1$ between internal groupoids produces a dotted universal decomposition, as in the above right hand side diagram, where the internal functor $(u_0,u_1)$ is a discrete fibration.

In a category $\mathbb E$, any of these properties guarentees the existence of normalizers. They become a characterization of this existence when $\mathbb E$ is quasi-pointed and protomodular. The universal decomposition 1) in the category $Gp$ of groups is described in detail in Section \ref{groupe}, from which the decomposition 2) is straighforward with the end of Section \ref{xxxgp}. 

Examples of non-pointed categories which satisfy this property are given with any slice or coslice category of the pointed protomodular categories with normalizers, and, in various circumstances, with any fibre $Grd_Y\mathbb E$ of the fibration $(\;)_0:Grd\mathbb E\rightarrow \mathbb E$ of internal groupoids, for instance when $\mathbb E$ is a Mal'tsev category.

As a collateral effect, this new approach allows us to clarify the relationship between existence of normalizers and  Mal'tsevness or (strong) protomodularity of the ground category $\mathbb E$, see  Section \ref{structob}, Lemma \ref{malc1}, Corollaries \ref{malc2} and  \ref{strongprot}. It also sheds new light on the non-pointed additive setting (in the sense of \cite{B15}), giving rise to several subtle differentiations, see Section \ref{adset}.

The article is organized along the following lines:\\
Section 1) is devoted to introducing the notion of $\Theta$-extremal decomposition of monomorphisms which is our categorical conceptual setting leading to the existence of normalizers.
Section 2) investigates the particular case of the categories which have $\P$-extremal decomposition of monomorphisms (namely which satisfy the above universal decomposition 1)), and it determines their first properties: in particular they are showed to be necessarily Mal'tsev categories.
Section 3) is devoted to the proof of the equivalence between  the above universal decompositions 1) and 2).
Section 4) is devoted to the relationship between these universal decompositions and the existence of normalizers.
Section 5) is devoted to the stability of these decompositions under slicing and coslicing.
Section 6) is devoted to the characterization theorem associated with the quasi-pointed protomodular setting.
Finally section 7) is devoted to the relationship between some kinds of these decompositions and the  non-pointed additive setting.

All the results developped here appeared in the long preprint \cite{BSM} were many other results are given, among them the fact that, in the same way as the existence of normalizers in the pointed case \cite{BG2}, the $\P$-extremal decomposition implies that: 1) any equivalence relation $R$ has a centralizer (i.e \emph{action distinctiveness}) and: 2) in the protomodular context, any subobject in a fibre $Pt_Y\EE$ has a centralizer (i.e \emph{fiberwise algebraic cartesian closedness}). 

\section{$\Theta$-extremal decomposition}

In this article any category $\EE$ will be supposed finitely complete.
Let $\Theta$ be a class of morphisms in  $\mathbb E$: it is said to be \emph{quasi-proper} when it contains the isomorphisms, is stable under composition with them, and is stable under product and pullback; it is said to be \emph{proper} when it contains the isomorphisms, is stable under composition and pullback and is such that, whenever $g.f$ and $g$ are in $\Theta$, the map $f$ is in $\Theta$.  It is clear that a proper class is quasi-proper. 
\begin{defi}\label{def1}
Let $\Theta$ be any class in a category $\mathbb E$, and $v:U\into T$ a monomorphism. We say that a decomposition $v=w.u$ with $u\in \Theta$ is $\Theta$-extremal, when any other decomposition $v=w'.u'$ with $u'$ in $\Theta$:
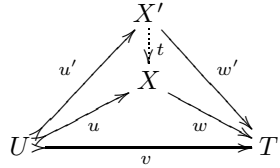
$$\xymatrix@=12pt
{  && X' \ar[ddrr]^{w'} \ar@{.>}[d]^{t} \\
   && X  \ar[drr]_{w}  \\
   {U\;} \ar@{>->}[rru]_{u} \ar@{>->}[rruu]^{u'} \ar@{>->}[rrrr]_{v} &&&& T
}
$$
determines a unique factorization $t$. A monomorphism $v:U\rightarrowtail T$ in $\mathbb E$ will be said to be a $\Theta$-outsider when $v=v.1_U$ is an extremal decomposition.
\end{defi}
\noindent Clearly the map $u$ is a monomorphism, and a monomorphism $v:U\rightarrowtail T$ is in $\Theta$ if and only if $v=1_T.v$ is an extremal decomposition. A map in $\Theta$ which is also a $\Theta$-outsider is an isomorphism. We shall now investigate the first properties of this  kind of decomposition.
\begin{lemma}\label{lemmaout}
	Suppose the monomorphism $v:U\rightarrowtail T$ is a $\Theta$-outsider. Then $v$ has no other monomorphic decomposition $v=w.u$ with $u$ in $\Theta$ than $v=v.1_U$, up to isomorphism.
\end{lemma}
\proof
Suppose $v=u'.w'$ with $w'$ monomorphic, then the factorization $\tau$:
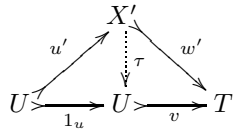
$$ \xymatrix@=20pt{
	& {X'\;} \ar@{>->}[rd]^{w'}\ar@{.>}[d]^{\tau} \\
	{U\;} \ar@{>->}[r]_{1_u} \ar@{>->}[ru]^{u'}& {U\;} \ar@{>->}[r]_v  & T 
}
$$
is a monomorphism since so is $w'$; then it is an isomorphism since it is split by $u'$ and $u'$ is an isomorphism as well.
\endproof
The first important observation is the following one:
\begin{lemma}\label{mon}
Suppose $\Theta$ quasi-proper. If $v=w.u$ is a $\Theta$-extremal decomposition, then the coreflector $w$ is necessarily a monomorphism. 
\end{lemma}
\proof
Complete the following right hand side square with the kernel equivalence relations $R[v]$ and $R[w]$ of the maps $v$ and $w$:
$$ \xymatrix@=20pt{
     R[v]=U \ar@<-1ex>[r]_>>>>{1_U} \ar@<1ex>[r]^>>>>{1_U} \ar[d]_{R(u)} & {U\;}   \ar[l]  \ar@{>->}[r]^v \ar[d]^u & T \ar@{=}[d]\\
     R[w] \ar@<-1ex>[r]_{p_{1}} \ar@<1ex>[r]^{p_{0}} & X   \ar[l]  \ar[r]_w  & T
                   }
$$
Since the left hand side part of the diagram is a joint pullback and $u$ is in $\Theta$ which is quasi-proper, then $R(u)$ is in $\Theta$. Accordingly the decomposition $v=(w.p_0).R(u)$ produces a unique factorization through $w$. So, we get $p_0=p_1$, and $w$ is a monomorphism.
\endproof
\begin{lemma}\label{suc}
Suppose $\Theta$ is quasi-proper. Let $v=v_2.v_1$ be a monomorphic decomposition of $v$ and $v=w.u$ a $\Theta$-extremal decomposition. Then the $\Theta$-extremal decomposition $v_1=w_1.u_1$  is given by the pullback of the $\Theta$-extremal decomposition $v=w.u$ along the monomorphism $v_2$.
\end{lemma} 
\proof
The map $v_2$ being a monomorphism, the following vertical rectangle is pullback. Introduce the lower quadrangle as a pullback and denote $u_1$ the canonical factorization:
$$ \xymatrix@=12pt{
   &  {U_1\;} \ar@{=}[rrr] \ar@{ >->}[dd]_{v_1} \ar@{>->}[dr]^{u_1} \ar@{>->}[dl]_{u'_1} &&& {U_1\;} \ar@{ >->}[dd]_<<<<{v_2.v_1} \ar@{>->}[dr]^{u}\\
   {X'_1\;} \ar@{>->}[dr]_{w'_1} & & {X_1\;} \ar@{>->}[dl]^{w_1} \ar@{>->}[rrr]_x &&&  {X\;} \ar@{>->}[dl]^{w} \\
   &  {U_2\;} \ar@{>->}[rrr]_{v_2} &&& T
                   }
$$
So, the upper quadrangle is a pullback and $u_1$ is a monomorphism in $\Theta$. Let $v_1=w'_1.u'_1$ be a decomposition with $u'_1$ a monomorphism in $\Theta$. Then the decomposition $v=v_2.v_1=(v_2.w'_1).u'_1$ produces a factorization $t:X'_1\rightarrow X$ satisfying  $w.t=v_2.w'_1$ which assures the factorization through the vertex $X_1$ of the lower quadrangled pullback. 
\endproof
The previous lemma gives immediately rise to the following one:
\begin{lemma}\label{suc1}
\noindent Suppose $\Theta$ is quasi-proper. Let $w=m.w'$ be any monomorphic decomposition and $w$ a $\Theta$-outsider. Then $w'$ is a $\Theta$-outsider.
\end{lemma}
Straighforward is the following one as well:
\begin{lemma}\label{outsider}
If $\Theta$ is is stable under composition and, a fortiori, if it is proper, when $v=w.u$ is a $\Theta$-extremal decomposition, the map $w$ is a $\Theta$-outsider.
\end{lemma}
\noindent Here is the main tool of this work: 
\begin{defi}\label{defi2}
Let $\Theta$ be a quasi-proper class in a category $\mathbb E$. This category will be said to have $\Theta$-extremal decompositions of monomorphisms when any monomorphism $v:U\rightarrowtail T$ has a $\Theta$-extremal decomposition $v=w.u$ with respect to the class $\Theta$. It will be said to have stably $\Theta$-extremal decompositions of monomorphisms when these $\Theta$-extremal decompositions are stable under pullbacks along maps in $\Theta$. 
\end{defi}
Warning: the $\Theta$-extremal decompositions do not preserve the inclusion of subobjects.

When $\mathbb E$ has stably $\Theta$-extremal decompositions of monomorphisms, the $\Theta$-outsider monomorphisms are stable under pullback along maps in $\Theta$. 

\smallskip

Define a normal monomorphism in the category $Gp$ of groups as an injective homomorphism $m:H\rightarrowtail G$ such that $m(H)$ is a normal subgroup of $G$. The class $N$ of normal monomorphisms is quasi-proper in $Gp$ but not proper.
\begin{ex}
The category $Gp$ has $N$-extremal decompositions of monomorphisms which are not stable.
\end{ex}
\proof
Starting from any monomorphism $n:H\rightarrowtail G$, its $N$-extremal decomposition is given by $\tilde n:H\stackrel{\simeq}{\rightarrowtail} n(H) \rightarrowtail N(n(H))$ where $N(G')$ is the normalizer of the subgroup $G'\rightarrowtail G$. 
\endproof
In the same way, define an ideal monomorphism in the category $Rg$ of non-unitary rings (resp. $K$-$Lie$ of Lie-algebras) as an injective homomorphism $m:B\rightarrowtail A$ such that $m(B)$ is an ideal of $A$. The class $I$ of ideal monomorphisms is quasi-proper in $Rg$ but not proper.
\begin{ex}
	The category $Rg$ (resp. $K$-Lie) has $I$-extremal decompositions of monomorphisms which are not stable.
\end{ex}

\begin{prop}\label{osuc}
Suppose $\Theta$ is proper and $\mathbb E$ has $\Theta$-extremal decompositions of monomorphisms. When, in addition, $\Theta$ satisfies  three out of two condition, any decomposition $v=w.u$ of the monomorphism $v$ with $u$ in $\Theta$ and $w$ a $\Theta$-outsider monomorphism is $\Theta$-extremal.
\end{prop}
\proof
Let $v=w'.u'$ be the $\Theta$-extremal decomposition:
$$ \xymatrix@=20pt{
    & {X'\;} \ar@{>->}[rd]^{w'} \\
    {U\;} \ar@{>->}[r]_u \ar@{>->}[ru]^{u'}& {X\;} \ar@{>->}[r]_w \ar@{.>}[u]^{\bar u} & T 
}
$$
There is a factorization $\bar u$ which is a monomorphism since so is $w$ and which, according to Lemma \ref{suc1}, is a $\Theta$-outsider since so is $w$.
Since, moreover, $\Theta$ satisfies the \emph{three out of two} condition (namely: if two maps among the triple $(u,\bar u,u'=\bar u.u)$ are in $\Theta$, so is the third one), then $\bar u$ is in $\Theta$; accordingly it is an isomorphism.
\endproof
We shall be specially interested by the following situation which will give us, in the next section, many examples of categories having \emph{stable} $\Theta$-extremal decompositions.
\begin{defi}\label{fibfib}
	Let $F: \mathbb D \to \mathbb C$ be a left exact fibration. We shall say that the ground category $\mathbb C$ is (resp. stably) $F$-decomposable when $\mathbb D$ has (resp. stably) extremal decompositions of monomorphisms with respect to the proper class of $F$-cartesian maps.
\end{defi}

\section{$\P$-decomposable categories}

Given a category $\mathbb E$, recall \cite{B-1} that $Pt\mathbb E$ denotes the category whose objects are the split epimorphisms (where split epimorphism means split epimorphism with a given splitting) in $\mathbb{E}$ and whose arrows are the commuting squares between such split epimorphisms, and that $\P_{\mathbb{E}} \colon Pt\mathbb E \to \mathbb{E}$ denotes the functor associating with any split epimorphism its codomain: it is \emph{the fibration of points}. The $\P_{\mathbb{E}}$-cartesian maps are nothing but the pullbacks of split epimorphisms and determine a proper class in $Pt\mathbb E$ we shall sometimes denote by $\P$ for short. As we shall recall below this fibration $\P_{\mathbb{E}}$ has strong classfication properties.

\subsection{Some structural observations}\label{structob}

Let us begin with a first clarification related to the $\P_\EE$-outsider monomorphims. For that, let us recall that a Mal'tsev category is a category in which any reflexive relation is an equivalence relation, see \cite{CLP} and \cite{CPP}, and that a protomodular category is a category such that any  base-change functor with respect to the fibration  $\P_{\mathbb{E}}$ is conservative. A category $\EE$ is protomodular if and only if the class of $\P$-cartesian morphisms satisfies the three out of two property \cite{B-1} (see definition in the proof of Proposition \ref{osuc}).  Any protomodular category is a Mal'tsev one. 
\begin{prop}\label{malc1}
	Any  category $\mathbb E$ in which any $\P_{\EE}$-invertible monomorphism is a $\P_{\EE}$-outsider is necessarily a Mal'tsev category. In a protomodular category $\EE$, any $\P_{\EE}$-invertible monomorphism is necessarily a $\P_{\EE}$-outsider, and any decomposition in $Pt\EE$ with a left hand side pullback is necessarily $\P_{\EE}$-extremal: 
	$$ \xymatrix@=8pt{
		&  {X'\;} \ar@{>->}[rrr]^{x} \ar@<-1ex>[dd]_{f'}  &&& {\bar X\;} \ar@<-1ex>[dd]_{\bar f} \ar@{>->}[rrr]^{\bar x} &&& X \ar@<-1ex>[dd]_{f} \\
		& &  &&&  \\
		&  {Y'\;} \ar@{>->}[rrr]_{y}  \ar[uu]_{s'} &&& Y  \ar[uu]_{\bar s} \ar@{=}[rrr] &&& Y \ar[uu]_{s}
	}
	$$
\end{prop}
\proof
Let be given  a reflexive relation $(d_0,d_1):R\rightrightarrows X$ on the object $X$. Then consider the following left hand side commutative diagram in $Pt\mathbb E$ where, by assumption, the whole rectangle is a $\P$-outsider:
$$
\xymatrix@=20pt{
	{R\;}\ar@{>->}[r]^{s_1} \ar@<-1ex>[d]_{d_0}\ar@{ >->}@(u,u)[rrr]^{(d_0,d_1)} & {R[d_0]\;} \ar@<-1ex>[d]_{p_0} \ar[rr]^{(d_1.p_0,d_1.p_1)} &&  X\times X \ar@<-1ex>[d]_{p_0}&& {R[d_0]\;} \ar@<-1ex>[d]_{d_0} \ar[r]^{d_2} &  R \ar@<-1ex>[d]_{d_0}\\
	{X\;} \ar@{>->}[r]_{s_0} \ar[u]_{s_0}\ar@{=}@(d,d)[rrr]_{} & {R\;} \ar[u]_{s_0} \ar[rr]_{d_1}  && X \ar[u]_{s_0} && {R\;} \ar[r]_{d_1} \ar[u]_{s_0} & {X\;} \ar[u]_{s_0} 
}             
$$
Since its left hand side part is a pullback, i.e. is a $\P$-cartesian monomorphism, there is a factorization $(d_1,d_2)$ given as on the right hand side satisfying $(d_0,d_1).d_2=(d_1.p_0,d_1.p_1)$; this shows $R[d_0]\subset (d_1)^{-1}(R)$, so that $R$ is an equivalence relation.

Now suppose that $\EE$ is protomodular and take any $\P_{\EE}$-invertible monomorphism $(1_Y,i):(f',s')\into (f,s)$. Now, given any decomposition $(1_Y,i)=(\beta,\bar{\beta}).(\alpha,\bar{\alpha}) (*)$ where $(\alpha,\bar{\alpha})$ is monomorphic and cartesian, consider the following pullback in $Pt\EE$:
$$\xymatrix@=5pt{
	P \ar[rrd]_{\check{\beta}}  \ar[ddd]_>>>>>{g'}   \ar@{ >->}[rrrrr]^{j} &&&&& W \ar[rrd]^{\bar{\beta}}  \ar[ddd]_>>>>>{g}\\
	&& X' \ar[ddd]_<<<<{f'}  \ar@{ >->}[rrrrr]^<<<<<{i}  &&&&& X \ar[ddd]_{f}  && \\
	&&& &&&&&&\\ 
	Z \ar[rrd]_{\beta}  \ar@{=}[rrrrr] \ar@<-1ex>[uuu]_<<<<<{t'}  &&&&& Z \ar@<-1ex>[uuu]_<<<<<{t} \ar[rrd]^{\beta}  \\
	&&  Y  \ar@{=}[rrrrr]  \ar@<-1ex>[uuu]_>>>>{s'}  &&&&& Y \ar@<-1ex>[uuu]_{s}
}
$$
The factorization $(*)$ in $Pt\EE$ produces the following factorization $(\alpha,\check{\alpha})$:
$$
\xymatrix@=20pt{
	{X'\;}\ar@{>->}[r]^{\check{\alpha}} \ar@<-1ex>[d]_{f'}\ar@{ >->}@(u,u)[rrr]^{\bar{\alpha}} & {P\;} \ar@<-1ex>[d]_{g'} \ar@{>->}[rr]^{j} &&  W \ar@<-1ex>[d]_{g}\\
	{Y\;} \ar@{>->}[r]_{\alpha} \ar[u]_{s'}\ar@(d,d)[rrr]_{\alpha} & {Z\;} \ar[u]_{t'} \ar@{=}[rr]  && Z \ar[u]_{t}  
}             
$$
Since  $(\alpha,\bar{\alpha})$ is cartesian (=underlying a pullback) and $j$ is a monomorphism, the left hand side square is a pullback as well. Now, when $\EE$ is protomodular, the three of of two condition for the $\P$-cartesian morphisms makes the right hand side diagram a  pullback and $j$ an isomorphism. With $(\beta,\check{\beta}.j^{-1}): (g,t) \to (f',s')$ we get the desired factorization which makes $(1_Y,i)$  a $\P_{\EE}$-outsider.

As for the last point, consider any other decomposition with $(v,u)$ $\P_{\mathbb E}$-cartesian:
$$ \xymatrix@=8pt{
	&& U \ar@<-1ex>[dd]_<<<{g} \ar[rrrrd]^{u'} &&& \\
	{X'\;} \ar@{>->}[rrr]_<<<<{x} \ar@<-1ex>[dd]_{f'} \ar@{>->}[rru]^{u} &&& {\bar X\;} \ar@<-1ex>[dd]_{} \ar@{>->}[rrr]_{\bar x} &&& X \ar@<-1ex>[dd]_{f} \\
	&& V \ar[uu]_>>>{t} \ar[rrrrd]^{v'} &&&  \\
	{Y'\;} \ar@{>->}[rrr]_{y}  \ar[uu]_{s'} \ar@{>->}[rru]^{v} &&& Y  \ar[uu]_{} \ar@{=}[rrr] &&& Y \ar[uu]_{s}
}
$$
Take the pullback $\bar u$ of $\bar x$ along $u'$; it determines a monomorphism in $Pt_{V}\mathbb E$ whose image by the change of base functor $v^*$ is the isomorphism $1_{X'}$ in $Pt_{Y'}\mathbb E$ since $(v,u)$ is $\P_{\mathbb E}$-cartesian. Now, $\mathbb E$ being protomodular, $\bar u$ is an isomorphism which produces the desired factorization.
\endproof

\begin{prop}\label{thinvert}
	Let $\mathbb E$ be a $\P$-decomposable category. If the monomorphism $(\gamma,\bar{\gamma}):(a',b')\into (a,b)$ in $Pt\mathbb E$ has a decomposition $(w',\bar w').(u',\bar u')$ where $(w',\bar w')$ is $\P_{\mathbb E}$-invertible and $(u',\bar u')$ is a $\P_{\mathbb E}$-cartesian monomorphism, then its $\P_{\EE}$-extremal decomposition $(w,\bar w).(u,\bar u)$ is such that $(w,\bar w)$ is $\P_{\mathbb E}$-invertible as well. Accordingly any $\P_{\mathbb E}$-invertible monomorphism in $Pt\mathbb E$ is a $\P_{\mathbb E}$-outsider.
\end{prop}
\proof
Consider the following diagrams of split epimorphisms where $w'$ is invertible and the right hand side diagram is an extremal decomposition in $Pt\mathbb E$:
$$
\xymatrix@=20pt{
	{\bar U\;} \ar[d]_{a'} \ar@(u,u)@{>->}[rr]^{\bar{\gamma}} \ar@{>->}[r]^{\bar u'} & {\bar A'\;} \ar[d]_{\bar a'} \ar[r]^{\bar w'}  & \bar T \ar[d]_{a} && {\bar U\;} \ar[d]_{a'} \ar@(u,u)@{>->}[rr]^{\bar{\gamma}} \ar@{>->}[r]^{\bar u} & {\bar A\;} \ar[d]_{\bar a} \ar@{>->}[r]^{\bar w}  & \bar T \ar[d]_{a}\\
	{U\;} \ar@<-1ex>[u]_{b'} \ar@(d,d)@{>->}[rr]_{\gamma} \ar@{>->}[r]_{u'} & {A'\;} \ar[r]_{w'}^{\simeq} \ar@<-1ex>[u]_{\bar b'} & T \ar@<-1ex>[u]_{b} && {U\;} \ar@<-1ex>[u]_{b'} \ar@(d,d)@{>->}[rr]_{\gamma} \ar@{>->}[r]_{u} & {A\;} \ar@{>->}[r]_{w} \ar@<-1ex>[u]_{\bar b} & T \ar@<-1ex>[u]_{b}
}
$$
Then there is a factorization:
$$
\xymatrix@=20pt{
	{\bar A'\;} \ar[d]_{\bar a'} \ar[r]^{\bar t} & {\bar A\;} \ar[d]_{\bar a}  \\
	{A\;} \ar@<-1ex>[u]_{\bar b'} \ar[r]_{t} & {A\;}  \ar@<-1ex>[u]_{\bar b} 
}
$$
such that $w.t=w'$. Now, since $w'$ is an isomorphism and $w$ is a monomorphism, the map $w$ is an isomorphism. The last assertion is then straightforward.
\endproof

Now, with Proposition \ref{malc1}, we get immediately the following:

\begin{coro}\label{malc2}
	Any $\P_{\mathbb E}$-decomposable category $\mathbb E$ is a Mal'tsev one. 
\end{coro}
 
 \subsection{Examples of $\P_{\EE}$-decomposable categories}\label{groupe}

As already said, Definition \ref{fibfib} will provide us with many examples of stably $\Theta$-decomposable categories: as a first step, let us show that the categories $Gp$, $Rg$ and $K$-$Lie$ are stably $\P$-decomposable. Let us begin by describing the extremal decomposition in $Pt(Gp)$. Given any subobject $(j,i)$ in this category:
$$
\xymatrix@=20pt{
	{X'\;}\ar[d]_{f'} \ar@(u,u)@{ >->}[rr]_{i} \ar@{>.>}[r]^{} & {\bar X\;} \ar[d]_{\bar f} \ar@{>.>}[r]^{}  & X \ar[d]_{f}\\
	{Y'\;}\ar@<-1ex>[u]_{s'} \ar@(d,d)@{ >->}[rr]^{j} \ar@{>.>}[r]_{} & {\bar Y\;} \ar@{>.>}[r]_{} \ar@<-1ex>[u]_{\bar s} & Y \ar@<-1ex>[u]_{s}
}
$$
the subgroup $\bar Y$ is $\{y\in Y/ s(y).u.s(y)^{-1}\in X', \;\forall u\in Kerf'\}$ while the subgroup $\bar X$ is $\{  \;{x\in X/ x\in f^{-1}(\bar Y)\; \rm and}\; x.sf(x)^{-1}\in X'\}$. \emph{Let us show that it is a $\P$-extremal decomposition}
\proof
The subset $\bar Y$ is clearly a subgroup of $Y$. And $\bar X$ is clearly stable under inversion. Now suppose $a$ and $b$ in $\bar X$. Then $a.b$ is in $f^{-1}(\bar Y)$. Moreover $b.sf(b)^{-1}$ is in $Ker f'$ and $f(a)$ in $\bar Y$, so that:\\ $a.b.sf(b)^{-1}.sf(a)^{-1}=(a.sf(a)^{-1}).(sf(a).b.sf(b)^{-1}.sf(a)^{-1})$ is in $Ker f'$.\\ So, $a.b$ in $\bar X$, and $\bar X$ is a subgroup of $X$.

The left hand side square is a pullback: suppose $a\in \bar X$ and such that $f(a)\in Y'$. Then we have $sf(a)\in X'$ and $a=(a.sf(a)^{-1}).sf(a)\in X'$.

Now consider a commutative diagram where the left hand side square is a pullback: 
$$
\xymatrix@=20pt{
	{X'\;}\ar[d]_{f'} \ar@(u,u)@{ >->}[rr]_{i} \ar@{>->}[r]^{} & {W\;} \ar[d]_{g} \ar[r]^{k}  & X \ar[d]_{f}\\
	{Y'\;}\ar@<-1ex>[u]_{s'} \ar@(d,d)@{ >->}[rr]^{j} \ar@{>->}[r]_{} & {Z\;} \ar[r]_{h} \ar@<-1ex>[u]_{t} & Y \ar@<-1ex>[u]_{s}
}
$$
First, let us show that $h(z)$ is in $\bar Y$. Let $u$ be in $Ker f'$, then we have to show that $sh(z).u.sh(z)^{-1}=kt(z).u.kt(z)^{-1}=k(t(z).u.t(z)^{-1})$ is in $X'$. It is enough to show that the image by $g$ of $t(z).u.t(z)^{-1}$ is in $Y'$. Now $g(t(z).u.t(z)^{-1})=z.1.z^{-1}=1$ which is in $X'$. It remains to show that $k(w)$ is in $\bar X$. 1) $fk(w)=hg(w)$ is in $\bar Y$ according to our first step. 2) $k(w).sfk(w)^{-1}=k(w.tg(w)^{-1})$ is in $X'$ as soon as $g(w.tg(w)^{-1})=1$ is in $Y'$ which is straightforward.
\endproof
 
Similarly: 1) in $Rg$ the extremal decomposition is obtained  in the following way: $\bar Y=\{y\in Y/ s(y).u\;  {\rm and} \; u.s(y)\in X', \;\forall u\in Kerf'\}$ while $\bar X=\{  \;{x\in X/ x\in f^{-1}(\bar Y)\; \rm and}\; x-sf(x)\in X'\}$;\\ 2) in $K$-$Lie$ it is obtained  in the following way: $\bar Y=\{y\in Y/ [s(y),u]
\in X', \;\forall u\in Kerf'\}$ while $\bar X=\{  \;{x\in X/ x\in f^{-1}(\bar Y)\; \rm and}\; x-sf(x)\in X'\}$.

These three previous categories are pointed. Now, the above description for the category $Rg$ remains valid for the category  $Rg_*$ of unitary rings which is longer a pointed one.

\emph{It remains to show that these four categories are stably $\P$-decomposable}. 
\proof Since any protomodular category satisfies the  three out of two condition and since the four categories are  protomodular, it is enough to show that the $\P$-ousider monomorphims are stable under pullback along the $\P$ cartesian maps, according to Proposition \ref{osuc}. It is straigforward to show that in the four cases, a monomorphism $(j,i)$ is a $\P$-outsider if and only if $\bar Y=Y'$, from which $\bar X=X'$ follows. 
Consider the following pullback in $Pt\EE$ where the front square is a pullback:
$$\xymatrix@=5pt{
	k^{-1}(X') \ar@{ >->}[rrd]_{}  \ar[ddd]_>>>>>{\check f'}   \ar[rrrrr]^{} &&&&& X' \ar@{ >->}[rrd]^{i}  \ar[ddd]_>>>>>{f'}\\
	&& \check X \ar[ddd]_<<<<{\check f}  \ar[rrrrr]^<<<<<{k}  &&&&& X \ar[ddd]_{f}  && \\
	&&& &&&&&&\\ 
	h^{-1}(Y') \ar@{ >->}[rrd]_{}  \ar[rrrrr] \ar@<-1ex>[uuu]_<<<<<{\check s'}  &&&&& Y' \ar@<-1ex>[uuu]_<<<<<{s'} \ar@{ >->}[rrd]^{j}  \\
	&&  \check Y  \ar[rrrrr]_h  \ar@<-1ex>[uuu]_>>>>{\check s}  &&&&& Y \ar@<-1ex>[uuu]_{s}
}
$$
The front square being a pullback, the restriction of the homomorphism $k$ to $Ker \check f$ produces an isomorphism $Ker \check f\simeq Ker f$ and an isomorphism  $Ker \check f'\simeq Ker f'$. We have to show that $\bar{\check Y}=h^{-1}(Y')$ whenever $\bar Y=Y'$.

In $Gp$: let $c\in\check Y$ be such that $\check s(c).v.\check s(c)^{-1}$ is in $k^{-1}(X')$ for  all $v\in Ker \check f'$. We have to show that $h(c)$ is in $Y'$. For that it is enough to show that $sh(c).u.sh(c)^{-1}$ is in $X'$ for all $u\in Ker f'$. According to our above remark, there is a unique $v\in Ker \check f'$ such that $u=k(v)$. So, $sh(c).u.sh(c)^{-1}=k(\check s(c).v.\check s(c)^{-1})$ with $\check s(c).v.\check s(c)^{-1}\in k^{-1}(X')$, whence $sh(c).u.sh(c)^{-1} \in X'$.

In $Rg$ and $Rg_*$: let $c\in\check Y$ be such that $\check s(c).v \; {\rm and} \; v.\check s(c)\; \in k^{-1}(X')$ for  all $v\in Ker \check f'$. We have to show that $h(c)$ is in $Y'$. For that it is enough to show that $sh(c).u\; {\rm and} \; u.sh(c)\in X'$ for  all $u\in Kerf'$. Again  there is a unique $v\in Ker \check f'$ such that $u=k(v)$. So, $sh(c).u=k(\check s(c).v)$ with $\check s(c).v\in k^{-1}(X')$, whence $sh(c).u \in X'$; the same proof holds for $u.sh(c)$.

In $K$-Lie: let $c\in\check Y$ be such that $[\check s(c),v] \in k^{-1}(X')$ for  all $v\in Ker \check f'$. We have to show that $h(c)$ is in $Y'$. For that it is enough to show that $[h\check s(c),u] \in X'$ for  all $u\in Ker f'$. Again  there is a unique $v\in Ker \check f'$ such that $u=k(v)$. So, $[sh(c),u]=k[\check s(c),v]$ with $[\check s(c),v]\in k^{-1}(X')$, whence $[sh(c),u] \in X'$.
\endproof

\subsection{The case of monoids and semi-rings}

The category $Mon$ of  monoids is not a Mal'tsev one; so, according to Corollary \ref{malc2}, it cannot be $\P_{\EE}$-decomposable. However it is $F$-decomposable for some subfibration $F$ of $\P_{\EE}$.

 In \cite{MMS} a split epimorphism $(f,s):X\rightleftarrows Y$ in $Mon$ is called a \emph{Schreier} split epimorphism when, for all $y\in Y$, the application  $\mu_y:Ker f\to f^{-1}(y)$ defined by $\mu_y(k)=k\cdot s(y)$ is bijective. This defines class $\Sigma$ of split epimorphims which is stable under pullback, and then determines a subfibration $\P_{\Sigma}$ of $\P_{\EE}$. Actually a split epimorphism $(f,s):X\rightleftarrows Y$ is a Schreier one if and only if there is a function $q: X\to Kerf$ such that $x=q(x).sf(x), \, \forall x\in X$ and $q(k.s(t))=k,\, \forall (k,t)\in Kerf \times Y$.
 
 \begin{prop}
 	The category $Mon$ is stably $\P_{\Sigma}$-decomposable.
 \end{prop}
 \proof
 Given any subobject in $Pt(Mon)$ between Schreier split epimorphims:
 $$
 \xymatrix@=20pt{
 	{X'\;}\ar[d]_{f'} \ar@<2ex>@{ >->}[rr]_>>>>>>{} \ar@{>.>}[r]^{} & {\bar X\;} \ar[d]_{\bar f} \ar@{>.>}[r]^{}  & X \ar[d]_{f}\\
 	{Y'\;}\ar@<-1ex>[u]_{s'} \ar@<-2ex>@{ >->}[rr]_{} \ar@{>.>}[r]_{} & {\bar Y\;} \ar@{>.>}[r]_{} \ar@<-1ex>[u]_{\bar s} & Y \ar@<-1ex>[u]_{s}
 }
 $$
 define the submonoid $\bar Y$ by $\{y\in Y/ q(s(y).u)\in X', \;\forall u\in Kerf'\}$ and the submonoid $\bar X$ by $\{x\in X/ x\in f^{-1}(\bar Y)\; {\rm and}\; q(x)\in X'\}$. The left hand side square is a pullback, since $x\in f^{-1}(Y')$ and $q(x)\in X'$ implies $x\in X'$. Whence $Ker \bar f=Ker f'$ where $(\bar f,\bar s)$ is the induced  split epimorphism. It is a Schreier one by taking the restriction of $q: X\to Kerf$ to $\bar X$ whose values are in $Ker f'(=Ker \bar f)$ by its definition.
 
 From this construction, checking its the universal property of this construction and its stablitily under pullback along $\P_{\Sigma}$-cartesian morphisms is straightforward.
 \endproof

\smallskip

Let $SRg$ be the category of semi-rings and $U:SRg\to CoM$ (where $CoM$ is the category of commutative monoids) the forgetful functor; it is left exact and conservative.  In \cite{MMS} a split epimorphism in $\bar{\Sigma}=U^{-1}(\Sigma)$ is, again called a Schreier one. Whence immediately:

 \begin{prop}
	The category $SRg$ is stably $\P{\bar\Sigma}$-decomposable.
\end{prop}

\section{Internal groupoids}\label{xxxgp}

In this section, we shall show that a category $\EE$ is $\P_{\mathbb E}$-decomposable if and only if the category $Grd\mathbb E$ of internal groupoids in $\mathbb E$ has $DiF$-decomposable, where $DiF$ is the class of discrete fibrations in $Grd\mathbb E$.

Let us recall that an internal groupoid $\underline Y_1$ is a reflexive graph as on the right hand side:
$$\xymatrix@=18pt{ \underline Y_1: & 
	R_2[d_0^Y] \ar[rr]^{p_{2}} \ar@<-2ex>[rr]|{p_{1}} \ar@<-4ex>[rr]_{p_{0}} \ar@(u,u)[rr]^{R(d_2^Y)}  && R[d_0^Y] \ar@(u,u)[rr]^{d_2^Y} \ar[rr]|{p_{1}} \ar@<-4ex>[rr]_{p_{0}} & & Y_1 \ar@<2ex>[ll]|{s_0} \ar@<-2ex>[ll]|{s_1} \ar@<2ex>[rr]^{d_1^Y} \ar@<-2ex>[rr]_{d_0^Y}  && Y_0 \ar[ll]|{s_0}
}
$$
endowed with a map $d_2^Y: R[d_0^Y]\rightarrow Y_1$ making the above  diagram a 3-truncated simplicial object.
In the set-theoretical context, we have $d_2^Y(\phi,\psi)=\psi.\phi^{-1}$. An internal functor is a morphism of 3-truncated simplicial objects. It is a discrete fibration (also called fibrant morphism) when any of the following leftward square indexed by $0$ (or equivalently $1$) is a pullback:
$$\xymatrix@=3pt{
 &	X_1 \ar[ddd]_{f_1}  \ar@<-1ex>[rrrr]_{d_{1}}\ar@<1ex>[rrrr]^{d_{0}} &&&&  {X_0\;} \ar[llll] \ar[ddd]^{f_0}  \\
	&  &&&&   & & &  &  && \\
	&&& &&&&&&\\
 &	Y_1 \ar@<-1ex>[rrrr]_{d_{1}} \ar@<1ex>[rrrr]^{d_{0}}    &&&&  {Y_0\;} \ar[llll]
}
$$
The class $DiF$ of fibrant morphisms is a proper class in $Grd\EE$.
\begin{theo}\label{lizer2}
The category $\mathbb E$ is (resp. stably) $\P$-decomposable if and only if the category $Grd\mathbb E$ of internal groupoids in $\mathbb E$ has (resp. stably) $DiF$-extremal decompositions of monomorphisms.
\end{theo}
\proof
The previous description of internal groupoids comes from the fact that $Grd\EE$ is actually the category of $T$-algebra, where $T$ is the following monad on the category $Pt\EE$, see \cite{B-2}:
$$
\xymatrix@=20pt{
	{X\;}\ar@{>->}[r]^{s_1} \ar@<-1ex>[d]_{f} & R[f] \ar@<-1ex>[d]_{d_0^f} & R^2[f] \ar@<-1ex>[d]_{d_0^f} \ar[l]_{d_2^f}\\
	{Y\;} \ar@{>->}[r]_{s}\ar[u]_{s}  & X \ar[u]_{s_0^f} & R[f] \ar[u]_{s_0^f} \ar[l]^{d_1^f}\\
	{(f,s)\;} \ar@{>->}[r]_{\lambda_{(f,s)}} & T(f,s)  & T^2(f,s)  \ar[l]^{\mu_{(f,s)}}	
}             
$$
So, the theorem will be the consequence of the following more general result.
\endproof
\begin{prop}
Let $U:\mathbb E\rightarrow \mathbb F$ be a left exact functor, $\Theta$ a proper class in $\mathbb F$. Then $\Theta'=U^{-1}(\Theta)$ is a proper class in $\mathbb E$. Suppose that $U$ has a left exact left adjoint $G$ such that $G.U$ preserves the maps in $\Theta$. Suppose moreover the natural transformation $\eta:1_{\mathbb F}\Rightarrow U.G$ is in $\Theta$. If the category $\mathbb E$ has (resp. stably) $\Theta'$-extremal decompositions of monomomorphisms, then the category $\mathbb F$ has (resp. stably) $\Theta$-extremal decompositions of monomorphisms as well.\\
If, moreover, the functor $U$ is monadic, the converse is true, namely if the category $\mathbb F$ has (resp. stably) $\Theta$-extremal decompositions on monomorphisms, then the category $\mathbb E$ has (resp. stably) $\Theta'$-extremal decompositions on monomorphisms as well. Moreover the functor $U$ preserves and reflects the extremal decompositions.
\end{prop}
\proof
The fact that $\Theta'$ is proper as soon as $\Theta$ is proper is straightforward. Now let $v:S\rightarrowtail T$ be a monomorphism in $\mathbb F$. Since $G$ is left exact, the map $G(v)$ is a monomorphism in $\mathbb E$. Let $m.n:G(S)\rightarrowtail W \rightarrowtail G(T)$ be its extremal decomposition in the category $\mathbb E$, the right hand side square below be a pullback in $\mathbb F$ and $u$ the induced factorization:
$$
\xymatrix@=20pt{
 {S\;} \ar[d]_{\eta_S} \ar@{>->}[r]_{u} \ar@(u,u)@{>->}[rr]^{v} &  {X\;} \ar[d]_{l} \ar@{>->}[r]_{w} & {T\;} \ar[d]^{\eta_T}   \\
 {U.G(S)\;} \ar@{>->}[r]^{U(n)} \ar@(d,d)@{>->}[rr]_{U.G(v)} & {U(W)\;}  \ar@{>->}[r]^{U(m)} & {U.G(T)\;}  
  }
$$
The map $\eta_T$ being in $\Theta$, so is $l$; the maps $\eta_S$ and $U(n)$ being in $\Theta$ (since $n$ is in $\Theta'=U^{-1}(\Theta)$), so is $u$. 

Let us show that $w.u$ is extremal in $\mathbb F$. Let $v=w'.u'$ with $u'$ a monomorphism in $\Theta$. The map $U.G(u')$ being in $\Theta$, the map $G(u')$ is in $\Theta'$ and the decomposition $G(v)=G(u').G(w')$ produces a factorization $t:G(X')\rightarrow W$ in $\mathbb E$ such that we have $m.t=G(w')$ and which, by adjonction, determimes a map $\tau: X'\rightarrow U(W)$ such that $U(m).\tau=\eta_T.w'$; whence the desired factorization $\bar{\tau}:X'\to X$. Suppose moreover $\mathbb E$ has stably $\Theta'$-extremal decompositions.  Starting with a map $\theta:\bar T \rightarrow T$ in $\Theta$, the pullback along $\theta$ in $\mathbb F$ preserves the previous construction, since $G(\theta)$ is in $\Theta'$ and $U$ left exact. Accordingly the category $\mathbb F$ has stably $\Theta$-extremal decompositions.

Conversely we shall show that when the monad $(T=U.G,\eta,\mu)$ on $\mathbb F$ is such that $T$ is left exact, preserves the maps in $\Theta$ and is such that $\eta:1_{\mathbb F}\Rightarrow T$ is in $\Theta$, the category $Alg^T$ has (resp. stably) $\Theta'$-extremal decompositions on monomorphisms  as soon as the category $\mathbb F$ has (resp. stably) $\Theta$-extremal decompositions on monomorphisms where $\Theta'$ is $(U^T)^{-1}(\Theta)$. So, let $v:(U,\alpha)\into (T,\beta)$ a monomorphism in $Alg^T$. Let us consider the following diagram, where the lower row is the extremal decomposition of $v$ in $\mathbb E$:
$$
\xymatrix@=20pt{
 {T(U)\;} \ar[d]_{\alpha} \ar@{>->}[r]_{T(u)} \ar@(u,u)@{>->}[rr]^{T(v)} &  {T(X)\;} \ar@{.>}[d]_{\xi} \ar@{>->}[r]_{T(w)} & {T(T)\;} \ar[d]^{\beta}   \\
 {U\;} \ar@{>->}[r]^{u} \ar@(d,d)@{>->}[rr]_{v} & {X\;}  \ar@{>->}[r]^{w} & {T\;}  
  }
$$
Let us show that the object $X$ of $\mathbb F$ is endowed with a $T$-algebra structure $\xi$. The monomorphisms $\eta_U$ and $T(u)$ are in $\Theta$. So that the decomposition $(\beta.T(w)).(T(u).\eta_U)=\beta.\eta_T.w.u=w.u=v$ produces a factorization $\xi:T(X)\to X$ which is easily seen to be a $T$-algebra structure since $w$ is monomorphic; this makes $v=w.u$ a decomposition in $Alg^T$; it is then straightforward to check that it is extremal for the class $\Theta'$. This construction shows that the functor $U$ preserves and reflects the extremal decompositions.\\
The stable aspect of this decomposition is straighforward from the left exactness of the endofunctor $T$. 
\endproof
So, $\mathbb E$ is $\P_{\mathbb E}$-decomposable if and only if any monomorphic internal functor $(v_0,v_1):\underline U_1\rightarrowtail \underline T_1$ between groupoids as on the left hand side diagram:
$$
\xymatrix@=20pt{
{U_1\;}\ar@{>->}[r]^{v_1} \ar@<1ex>[d]^{d_1}\ar@<-1ex>[d]_{d_0} & T_1 \ar@<1ex>[d]^{d_1}\ar@<-1ex>[d]_{d_0} && {U_1\;}\ar@{>->}[r]^{u_1} \ar@{>->}@(u,u)[rr]^{v_1} \ar@<1ex>[d]^{d_1}\ar@<-1ex>[d]_{d_0} & {X_1\;} \ar@<1ex>[d]^{d_1}\ar@<-1ex>[d]_{d_0} \ar@{>->}[r]^{w_1} &  T_1 \ar@<1ex>[d]^{d_1}\ar@<-1ex>[d]_{d_0}\\
{U_0\;} \ar@{>->}[r]_{v_0} \ar[u] & T_0 \ar[u] && {U_0\;} \ar@{>->}[r]_{u_0} \ar@{>->}@(d,d)[rr]_{v_0} \ar[u] & {X_0\;} \ar[u] \ar@{>->}[r]_{w_0}  & T_0 \ar[u]
}             
$$
produces an extremal decomposition as in the right hand side one, where the internal functor $(u_0,u_1)$ is a fibrant morphism. \emph{It is a kind of dual for the monomorphic functors of the comprehensive factorization for internal functors between groupoids} described in \cite{SW} and \cite{B-2}. Actually \emph{the result is even more precise}: if a monomorphism in $Pt\mathbb E$ is underlying a functor between groupoids as on the left hand side:
$$
\xymatrix@=20pt{
{U_1\;}\ar@{>->}[r]^{v_1} \ar@<1ex>@{.>}[d]^{d_1}\ar@<-1ex>[d]_{d_0} & T_1 \ar@<1ex>@{.>}[d]^{d_1}\ar@<-1ex>[d]_{d_0} && {U_1\;}\ar@{>->}[r]^{u_1} \ar@{>->}@(u,u)[rr]^{v_1} \ar@<1ex>@{.>}[d]^{d_1}\ar@<-1ex>[d]_{d_0} & {X_1\;} \ar@<1ex>@{-->}[d]^{d_1}\ar@<-1ex>[d]_{d_0} \ar@{>->}[r]^{w_1} &  T_1 \ar@<1ex>@{.>}[d]^{d_1}\ar@<-1ex>[d]_{d_0}\\
{U_0\;} \ar@{>->}[r]_{v_0} \ar[u] & T_0 \ar[u] && {U_0\;} \ar@{>->}[r]_{u_0} \ar@{>->}@(d,d)[rr]_{v_0} \ar[u] & {X_0\;} \ar[u] \ar@{>->}[r]_{w_0}  & T_0 \ar[u]
}             
$$
then the extremal decomposition in $Pt\mathbb E$ provides the middle vertical part with a unique factorization $d_1$ and a unique groupoid stucture, making fibrant the left hand side internal functor. Whence the immediate:
\begin{coro}
When $\EE$ is a $\P$-decomposable category, the following conditions are equivalent:\\
1) the monomorphic functor $(v_0,v_1)$ is a DiF-ousider;\\
2) its underlying monomorphism in $Pt\EE$ is a $\P$-outsider.\\
Accordingly any $(\;)_0$-invertible functor is necessarily a DiF-outsider.
\end{coro}

\section{Abstract normalizers}\label{abnorm}

In this section, we show in a very simple way how the above DiF-decomposition is related to the existence of normalizers. Let us recall the following:
\begin{defi}\label{norm}
A monomorphism $u$ in $\mathbb E$ is said to be normal to an equivalence relation $R$ when:\\
i) we have: $u^{-1}(R)=\nabla _U$\\
ii) the induced internal functor:
$$
\xymatrix@=20pt{
{\nabla_U=U\times U\;}\ar@{>->}[rr]^{\tilde u} \ar@<1ex>[d]^{p_1}\ar@<-1ex>[d]_{p_0} && R \ar@<1ex>[d]^{d_1}\ar@<-1ex>[d]_{d_0}\\
{U\;} \ar@{>->}[rr]_{u} \ar[u] && X \ar[u]
}             
$$
is a fibrant morphism in $Grd\EE$.
\end{defi}
In the category $Set$ of sets, when $U$ is not empty, it is equivalent to saying that $U$ is an equivalence class of $R$. Clearly, in this category, a monomorphism can be normal to many equivalence relations. In particular the inclusion $\emptyset\into X$ is normal to any equivalence relation $R$ on $X$ and in particular to $\nabla X$. However, in a protomodular category a monomorphim is normal to at most one equivalence relation, so that, for a monomorphism, being normal becomes a property \cite{B-1}. Recall now the following definition from \cite{BG2}:
\begin{defi}\label{orm}
	Given any category $\EE$, a monomorphism $v:U\into T$ has a normalizer when there is a pair $(u,R_v)$ with $u: U\into X$ normal to $R_v$ and a factorization $w:X\to T$ such that $v=w.u$ which is universal with respect to this kind of specific decomposition of $v$.
A category $\mathbb E$ is said to  have normalizers when  any monomorphism has a normalizer.
\end{defi}
By the universal property of a normalizer, this equivalence relation $R_v$ is the largest equivalence relation $R$ on $X$ to which the monomorphism $u$ is normal.

\begin{lemma}
When a monomorphism $v:U\into T$ has a normalizer, the induced factorization $w$ is necessarily a monomorphism.
\end{lemma}
\proof
Consider the following left hand side diagram and the right hand side one where $R[w]$ is the kernel equivalence relation of the map $w$:
$$
\xymatrix@=20pt{
	{\nabla_U\;} \ar@<1ex>[d]^{d_1^U} \ar@<-1ex>[d]_{d_0^U} \ar@{>->}[r]^{\tilde u} & {R_v\;} \ar@<1ex>[d]^{d_1^{R_v}} \ar@<-1ex>[d]_{d_0^{R_v}} \ar@{>->}[rr]^{s_0}  && R_v\times R_v \ar@<1ex>[d]^{d_1\times d_1} \ar@<-1ex>[d]_{d_0\times d_0} & {\nabla_U\;} \ar@<1ex>[d]^{d_1^U} \ar@<-1ex>[d]_{d_0^U} \ar@{>.>}[r]_{\sigma} \ar@(u,u)@{>->}[rr]^{s_0.\tilde u}& {\Sigma\;} \ar@<1ex>[d]^{d_1} \ar@<-1ex>[d]_{d_0} \ar@{>->}[r]_{(\delta_0^R,\delta_1^R)}  & R_v\times R_v \ar@<1ex>[d]^{d_1\times d_1} \ar@<-1ex>[d]_{d_0\times d_0}\\
	{U\;} \ar[u]  \ar@{>->}[r]_{u} & {X\;}\ar[u] \ar@{>->}[rr]_{s_0^X}  && X\times X  \ar[u]_{} & {U\;} \ar[u]_{}  \ar@{>->}[r]_{(u,u)} & {R[w]\;} \ar@{>->}[r]_{(d_0^w,d_1^w)} \ar[u]_{} & X\times X \ar[u]_{}
}
$$
The two lower maps are the same. So, if the right hand side part of the right hand side diagram is cartesian with respect to the fibration $(\;)_0:Grd\EE\to \EE$ (namely if $\Sigma$ is the inverse equivalence relation of $R_v\times R_v$ along $(d_0^w,d_1^w)$), there is a factorization $\sigma$ which induces a morphism $\nabla_U\into \Sigma$ of equivalence relation. We are going to show that this morphism is fibrant which will mean that the monomorphism $(u,u)$ is normal to $\Sigma$. First, it is cartesian since its composition with its right hand side cartesian neighbour is cartesian, this composition being equal to the whole left hand side diagram which is itself the composition of two cartesian maps. Now, since $w.u=v$ is a monomorphism, we get $u^{-1}(R[w])=\Delta_U$. Then consider the following left hand side commutative  diagram in $Grd\EE$ whose image by the fibration $(\;)_0$ is the right hand side square in $\EE$ which is a pullback by $u^{-1}(R[w])=\Delta_U$:
$$
\xymatrix@=20pt{
	{\nabla_U\;} \ar[d]_{s_0}  \ar@{>->}[rr]^{((u,u),\sigma)} && {\Sigma\;} \ar[d]^{}     && {U\;} \ar[d]_{s_0^U}  \ar@{>->}[r]^{(u,u)} & {R[w]\;} \ar[d]^{(d_0^w,d_1^w)}\\
	{\nabla_U\times \nabla_U\;}   \ar@{>->}[rr]_{(u,\tilde u)\times (u,\tilde u)} && {R_v\times R_v\;} && {U\times U\;}  \ar@{>->}[r]_{u\times u} & {X\times X\;} 
}
$$
Since the fibration $(\;)_0$ is left exact, the parallel horizontal arrows are cartesian and the image of this square by $(\;)_0$ is a pullback in $\EE$, then the left hand side square is itself a pullback in $Grd\EE$. So, since the lower functor is fibrant as the product of two fibrant morphisms, the upper one is fibrant as well.

Now we got a pair $((u,u),\Sigma)$ with $(u,u)$ normal to $\Sigma$ and two decompositions $v=w.(d_i^w.(u,u))$ ($i\in \{0,1\}$). Then, the universal property of the normalizer implies that $d_0^w=d_1^w$, and that, consequently $w$ is a monomorphism.
\endproof

\begin{prop}\label{surnorm}
	In a protomodular category $\EE$ any normal monomorphism $u:U\into X$ is its  own normalizer.
\end{prop}
\proof
Let $v$ be any normal monomorphism and $R_v$ the (unique) equivalence relation to which $v$ is normal in the protomodular category $\EE$. Let $v=g.u$ be any decomposition of $v$ with $u$ a normal monomorphism (to $R_u$). We have to show that $R_u$ can be factorized through $R_v$, which is equivalent to $R_u\subset g^{-1}(R_v)$. Clearly $u^{-1}(g^{-1}(R_v))=v^{-1}(R_v)=\nabla Y$. Then $R_u \cap g^{-1}(R_v)$ is normal to $u$ since it is included in $R_u$ and such that $u^{-1}(R_u \cap g^{-1}(R_v))=\nabla Y$. Now, since $\EE$ is protomodular, we get $R_u \cap g^{-1}(R_v)\simeq R_u$ and consequently $R_u \subset g^{-1}(R_v)$.
\endproof
\begin{prop}\label{lizer1}
Suppose $\mathbb E$ is $\P$-decomposable. Then any monomorphism $v:U\rightarrowtail T$ in $\mathbb E$ has a normalizer in the previous sense. Any normal monomorphism $u:U\into X$ in $\mathbb E$ admits a largest equivalence $R_u$ on $X$ to which $u$ is normal.
\end{prop}
\proof
It is a straightforward consequence of Theorem \ref{lizer2} and Proposition \ref{thinvert} applied to the following monomorphic morphism of equivalence relation:
$$
\xymatrix@=20pt{
  {U\times U\;} \ar@<-1ex>[d]_{p_0^U} \ar@<1ex>[d]^{p_1^U}\ar@<2ex>@{>->}[rr]^{v\times v} \ar@{>.>}[r]_{\tilde u} & {R\;} \ar@{>.>}[r]_{\tilde w}\ar@<-1ex>[d]_{d_0^R} \ar@<1ex>[d]^{d_1^R} & T\times T \ar@<-1ex>[d]_{p_0^T} \ar@<1ex>[d]^{p_1^T}\\
  {U\;} \ar[u] \ar@<-2ex>@{>->}[rr]_{v} \ar@{>.>}[r]^{u} & {X\;} \ar@{>.>}[r]^{w} & T \ar[u]
  }
$$
\endproof
More generally, let $R$ be any equivalence relation on the domain $U$ of a monomophism $v:U\into T$. Again, we can produce the monomorphism of equivalence relations given by the left hand side diagram:
$$
\xymatrix@=20pt{
	{R\;} \ar@<-1ex>[d]_{d_0^R} \ar@<1ex>[d]^{d_1^R}\ar@{>->}[rr]^{v\times v.(d_0^R,d_1^R)}    && T\times T \ar@<-1ex>[d]_{p_0^T} \ar@<1ex>[d]^{p_1^T} & {R\;}\ar@{>->}[r]^{\tilde u} \ar@{>->}@(u,u)[rr]^{v\times v.(d_0^R,d_1^R)} \ar@<1ex>[d]^{d_1^R}\ar@<-1ex>[d]_{d_0^R} & {S\;} \ar@<1ex>[d]^{d_1^S}\ar@<-1ex>[d]_{d_0^S} \ar@{>->}[r]^{\tilde w} &  T\times T \ar@<1ex>[d]^{p_1}\ar@<-1ex>[d]_{p_0} & {R\;}\ar@{>->}[rr]^{\tilde u}  \ar@<1ex>[d]^{d_1^R}\ar@<-1ex>[d]_{d_0^R} && {\bar R\;} \ar@<1ex>[d]^{d_1^{\bar R}}\ar@<-1ex>[d]_{d_0^{\bar R}}  \\
	{U\;} \ar[u] \ar@{>->}[rr]_{v}  && T \ar[u] & {U\;} \ar@{>->}[r]_{u} \ar@{>->}@(d,d)[rr]_{v} \ar[u] & {X\;} \ar[u] \ar@{>->}[r]_{w}  & T \ar[u] & {U\;} \ar@{>->}[rr]_{v}  \ar[u] && {T\;} \ar[u]
}
$$
When $\EE$ is $\P_\EE$-decomposable, the fibrant extremal decomposition given by the previous theorem,  produces the middle diagram, where $S$ is an equivalence relation on $X$ and the left hand side functor $(u,\tilde u)$ a fibrant monomorphism. This equivalence relation $S$ is the largest equivalence relation on $X$ which produces such a left hand side fibrant monomorphism. Clearly, this kind of decomposition can be extended  to any morphism of equivalence relations as on right hand side.

\smallskip

\noindent\textbf{The case of monoids and semirings}

\smallskip

Let us recall from \cite{MMS} the following:

\begin{defi}
	Given any category $\EE$ and any class $\Sigma$ of split epimorphims, we call $\Sigma$-equivalence relation any equivalence relation: $$\xymatrix@=7pt{
	 R  \ar@<2ex>[rr]^{d_1^R} \ar@<-2ex>[rr]_{d_0^R}  && X \ar[ll]|{s_0^R}
	}
	$$ such that the split epimorphism $(d_0^R,s_0^R)$ is in $\Sigma$. A morphism $f: X\to Y$ is called $\Sigma$-special when its kernel equivalence relation $R[f]$ is in $\Sigma$; an object $X$ is called $\Sigma$-special when the terminal map $X\to 1$ is $\Sigma$-special.
\end{defi}

Warning: a split $\Sigma$-special epimorphism is stronger than a split epimorphism in $\Sigma$; however an equivalence relation is a $\Sigma$-one, if and only if $(d_0^R,s_0^R)$ is a split $\Sigma$-special epimorphism. So if we denote by $\check{\Sigma}$ the class of split $\Sigma$-special epimorphisms, the monad $T$ of groupoids is stable on the subcategory $\P_{\check{\Sigma}}$. 

In $Mon$ and $SRg$ the Schreier-special objects are respectively the groups and the rings. According to the stability of the monad of groupoids on $\P_{\check{\Sigma}}$, in both cases, any morphism of Schreier equivalence relation  as on the left hand side produces a fibrant extremal decomposition of Schreier equivalence relations as on the right and  side:
$$\xymatrix@=20pt{
	{R\;} \ar@<-1ex>[d]_{d_0^R} \ar@<1ex>[d]^{d_1^R}\ar@{>->}[rr]^{\bar v}    && \bar R \ar@<-1ex>[d]_{d_0^{\bar R}} \ar@<1ex>[d]^{d_1^{\bar R}} && {R\;}\ar@{>->}[r]^{\tilde u} \ar@{>->}@(u,u)[rr]^{\bar v} \ar@<1ex>[d]^{d_1^R}\ar@<-1ex>[d]_{d_0^R} & {S\;} \ar@<1ex>[d]^{d_1^S}\ar@<-1ex>[d]_{d_0^S} \ar@{>->}[r]^{\tilde w} & \bar R \ar@<-1ex>[d]_{d_0^{\bar R}} \ar@<1ex>[d]^{d_1^{\bar R}}  \\
	{U\;} \ar[u] \ar@{>->}[rr]_{v}  && T \ar[u] && {U\;} \ar@{>->}[r]_{u} \ar@{>->}@(d,d)[rr]_{v} \ar[u] & {X\;} \ar[u] \ar@{>->}[r]_{w}  & T \ar[u] 
}
$$

\section{Slicing and coslicing}\label{slice}

In this section, we shall show that the $\P_\EE$-decomposition property is stable under slicing and coslicing. Let us recall that, given a category $\mathbb C$ and any object $Y$ in $\EE$, the slice category $\EE/Y$ is the category whose objects are the maps with codomain $Y$ and whose maps are the commutative triangles above $Y$. The coslice category $Y/\EE$ is defined by duality. The domain functor $dom:\EE/Y\rightarrow \mathbb \EE$ is discrete fibration which preserves and reflects pullbacks and equalizers while the codomain functor $cod:Y/\EE\rightarrow \EE$ is a discrete cofibration which preserves and relects pullbacks and equalizers as well.
\begin{prop}
	Let $U:\mathbb E\rightarrow \mathbb F$ be a  discrete fibration (resp. cofibration) which preserves and reflects pullbacks. Set $\Theta'=U^{-1}(\Theta)$. When $\mathbb F$ has (resp. stably) $\Theta$-extremal decompositions, then $\mathbb E$ has (resp. stably) $\Theta'$-extremal decompositions. When $\FF$ is (stably) $\FF$-decomposable, then $\EE$ is (stably) $\EE$-decomposable.
\end{prop}
\proof
Given any monomorphism $v:Z\rightarrowtail T$ in $\mathbb E$, consider the extremal decomposition $U(v)=w.u$ in $\mathbb F$. When $U$ is a discrete fibration (resp. cofibration) it determines a unique monomorphic decomposition $v=\bar w.\bar u$ above it since $U$ reflects the monomorphisms. The map $\bar u$ is in $\Theta'$ since $u$ is in $\Theta$. It is then straightforward that this decomposition is extremal. The fact that $U$ preserves and reflects pullbacks and the fact that we have $\Theta'=U^{-1}(\Theta)$, induces the assertion about the stability of this decomposition. For the last assertion, apply the first one to the functor $PtU:Pt\EE\to Pt\FF$.
\endproof
\begin{coro}\label{slic}
	The  $\P_\EE$-decomposable and stably $\P_\EE$-decomposable categories are stable under slicing and coslicing; accordingly they are both stable under the passage to any fibre $Pt_Y\mathbb E$, since $Pt_Y\mathbb E=1_Y/(\EE/Y)$. 
\end{coro}

So, the slice and coslice categories of the categories $Gp$ of groups, $R$ of rings and $R$-$Lie$ of Lie algebras on the ring $R$  produce new examples of non-pointed stably $\P$-decomposable categories. Starting with a (resp. stably) $\P$-decomposable category, any fiber $Pt_Y\EE$ becomes a pointed (resp. stably) $\P$-decomposable category. So, according to Proposition \ref{key1}, this fiber is a pointed category with normalizers where, accordingly, all the results of \cite{BG2} are valid.
\begin{prop}\label{presextr}
		If $\EE$ is stably $\P_\EE$-decomposable, the base-change functors with respect to the fibration $\P_{\EE}$ preserve the $\P$-extremal decompositions in the fibres.
\end{prop}
\proof
This comes from the two following observations:\\
1) given an object $Y$ in $\EE$ the domain functor $dom:Pt_Y\EE \rightarrow \EE$ reflects the extremal $\P_\EE$-decompositions;\\
2) given a map $h:Y'\rightarrow Y$ in $\EE$, the base-change functor $h^*:Pt_Y\EE \rightarrow Pt_{Y'}\EE$ produces $\P_{\EE}$-cartesian maps in $\EE$.\\
The conclusion follows from the assumed stable aspect of the $\P_\EE$-decomposition in $\EE$.
\endproof
\begin{coro}\label{strongprot}
	Let $\EE$ be stably $\P_\EE$-decomposable. If $\EE$ is protomodular, it is strongly protomodular.
\end{coro}
\proof
A protomodular category $\EE$ is strongly protomodular when, in addition, any (conservative) base-change functor $h^*:Pt_Y\EE \rightarrow Pt_{Y'}\EE$ reflects the normal monomorphisms; the categories $Gp$ of groups and $Rg$ of rings are examples of such categories, see \cite{BB}. This corollary is a consequence of following lemma and of the previous proposition, when $\EE$ is stably $\P_\EE$-decomposable.
\endproof
\begin{lemma}
	Let $H:\mathbb E \to \mathbb F$ be a conservative left exact functor between protomodular categories which are $\P$-decomposable. Then the functor $H$ reflects the normal monomorphisms as soon as $H$ preserves the $\P$-extremal decompositions.
\end{lemma}
\proof
We have to show that if the image by $H$ of the monomorphism $v:U\into T$ is normal, then $v$ is itself is normal. Let $v=w.u$ the decomposition through the normalizer $u$ of $v$. Since $H$ is left exact and preserves the $\P$-extremal decompositions and since $H(v)$ is normal, then $H(w)$ is an isomorphism. Accordingly, $H$ being conservative, the morphism $w$ is an isomorphism as well, and $v$ is normal.
\endproof

\section{Quasi-pointed categories}

In Section \ref{abnorm}, we showed that any $\P$-decomposable category has normalizers. Here we shall investigate a context in which the two conditions are equivalent, see Theorem \ref{lizer}.

\subsection{The general case}

 A category $\mathbb E$ is pointed when the terminal object is also initial; it is said quasi-pointed when it has an initial object $0$ and when, in addition, the map $0\rightarrow 1$ is a monomorphism, which implies that any initial map $\alpha_X:0\rightarrowtail X$ is a monomorphism. So, the fibre $Pt_0\mathbb E$ becomes a pointed full subcategory of $\mathbb E$.
 
Clearly the category $Set$ of sets in quasi-pointed.
Given any category $\mathbb E$, consider the fibration $(\;)_0:Grd\mathbb E\rightarrow \mathbb E$ associating with any internal groupoid $\underline Y_1$ its ``object of objects'' $Y_0$; the fibre $Grd_1\mathbb E$ above $1$ is nothing but the category $Gp\mathbb E$ of internal groups in $\mathbb E$ which is pointed; \emph{any other fibre $Grd_Y\mathbb E$ is quasi-pointed}, its initial object being $\Delta_Y$ the discrete equivalence relation on $Y$ and its terminal object being $\nabla_Y$ the indiscrete equivalence on $Y$. Furthermore any fibre $Grd_Y\mathbb E$ is protomodular \cite{B-1}. It will follow from our characterization theorem, that any  fiber $Grd_Y$ is actually a stably $\P$-decomposable category.

Suppose $\mathbb E$ is quasi pointed; we call \emph{kernel} of a map $f:X\rightarrow Y$ the upper horizontal arrow in the following left hand side pullback where $\alpha_Y$ is the initial map:
$$ \xymatrix@=20pt{
      {K[f]\;}    \ar@{>->}[r]^{k_f} \ar[d] & X \ar[d]^f && {EnX\;}    \ar@{>->}[r]^{\epsilon_X} \ar[d] & X \ar[d] && {EnX\;}    \ar@{>.>}[r]^{\bar{\epsilon}_X} \ar@{>->}@(u,u)[rr]^{\epsilon_X}\ar[d] & X\times X \ar[d]_{p_0^X}     \ar[r]^{p_1^X} \ar[d] & X \ar[d]\\
       {0\;}   \ar@{>->}[r]_{\alpha_Y}  & Y && {0\;}   \ar@{>->}[r]_{}  & 1 &&  {0\;}   \ar@{>->}[r]_{\alpha_X}  & X  \ar[r]_{} & 1 
                   }
$$
For the special case of the terminal map, we use the notations of the middle pullback and call this kernel $En X$ the \emph{endosome} of $X$, while we denote by $\bar{\epsilon}_X$ the unique factorization making $\bar{\epsilon}_X$ the kernel of $p_0^X$. Clearly the subobject $\epsilon_X$ is normal to $\nabla_X$. The middle pullback determines a left exact functor $En: \mathbb E \to Pt_0\mathbb E$ which is a right adjoint to the inclusion $Pt_0\mathbb E \hookrightarrow \mathbb E$.
\begin{lemma}\label{rose}
Let $\mathbb E$ be a quasi-pointed category. Then if a monomorphism $v:U\into X$ is normal to an equivalence relation $R$, so is $u.\epsilon_U: En U\into X$.
When $\mathbb E$ is protomodular, the converse is true. In this case a monomorphism $v: U\into T$ has a normalizer as soon as $v.\epsilon_U: En U\into X$ has one.
\end{lemma}
\proof
Since any $\epsilon_U$ is normal to $\nabla U$, if $u: U\into X$ is normal to $R$, so is $u.\epsilon_U$ in any category. Conversely when $u.\epsilon_U$ is normal to $R$, then $u^{-1}(u.\epsilon_U)=\epsilon_U$ is normal to $u^{-1}(R)$. If, in addition, $\mathbb E$ is protomodular, we get $u^{-1}(R)=\nabla_U$, and thanks to the three out of two conditions for the fibrant morphisms of equivalence relations, $u$ is normal (to $R$). 
Suppose now that $v.\epsilon_U=w.\tilde u$ is the factorization through the normalizer $(\tilde u,R_{\tilde u})$ of the monomorphism $v.\epsilon_U$. Then, since $\epsilon_U$ is normal, we get a factorization $u: U \into X$ where $X$ is the codomain of $\tilde u$ such that $w.u=v$ and $u.\epsilon_U=\tilde u$. According to the second equality, $u$ is normal to $R_{\tilde u}$. From that and the first assertion, $(u,R_{\tilde u})$ is necessarily the normalizer of $v$.
\endproof

When $\mathbb E$ is quasi-pointed, we shall denote by $Kt\mathbb E$ the category of split exact sequences, namely of split epimorphisms with a chosen kernel, and by $K: Kt\mathbb E\rightarrow Pt_0\mathbb E$ the functor associating with any split exact sequence the domain of its kernel map. Not only the functor $K$ is left exact, but \emph{it creates pullbacks and equalizers}. We shall denote by $J: \mathbb E \to Kt\mathbb E$ the functor associating with any $X$ the following exact secquence:
$$
\xymatrix@=15pt{
  {EnX\;}  \ar@{>->}[r]^{\bar{\epsilon}_X} & X\times X  \ar[r]_{p_0^X} & {\;X} \ar@<-1ex>@{>->}[l]_{s_0^X}
  }
$$
Clearly we have $K.J=En$.

On the other hand, it is straightforward that the forgetful functor $H:Kt\mathbb E\rightarrow Pt\mathbb E$ associating with any split exact sequence its underlying split epimorphim is a fully faithful and essentially surjective, namely that it determines a weak equivalence of categories, making the functor $\P_{\mathbb E}.H$ a fibration. 
Now, we need the following:
\begin{defi}
Given any functor $K:\mathbb E\to \mathbb F$, a map $f:X \to Y$ is pre-cartesian (resp. cartesian) with respect to $K$ whence it is universal among the maps in $\EE$ with codomain $Y$ whose image by $K$ is $K(f)$ (resp. whose image by $K$ factorizes through $K(f)$).\\
A left exact functor $K:\mathbb E\to \mathbb F$ is said to be \emph{pre-fibrant on monomorphisms} (resp. \emph{fibrant on monomorphisms}) when, given any monomorphism $u: U\into K(X)$ in $\FF$, there is a monomorphic pre-cartesian (resp. cartesian) map $\tilde u: \tilde U\into X$ in $\EE$ whose image by $K$ is isomorphic to $u$. 
\end{defi}
We have the straightforward following:
\begin{lemma}\label{key6}
Given any left exact functor $K:\mathbb E\to \mathbb F$, the following conditions are equivalent:\\
1) $K$ is conservative and pre-fibrant on monomorphisms;\\
2) $K$ is fibrant on monomorphisms;\\
3) $K$ determines a bijection between the set of isomorphic classes of subobjects of $K(X)$ and the set of isomorphic classes of subobjects of $X$.\\
Then $K$ determines a bijection between the set of isomorphic classes of equivalence relations on $K(X)$ and the set of isomorphic classes of equivalence relations on $X$.
\end{lemma}
\proof
This follows from the fact that any left exact functor is conservative as soon as it is conservative on monomorphisms, and that any left exact conservative functor is such that any monomorphism is cartesian. Then the point 3) determines a bijection between the set of isomorphic classes of reflexive relations on $K(X)$ and the set of isomorphic classes of reflexive relations on $X$. Now, $K$ being conservative, it reflects the equivalence relations among the reflexive ones since a reflexive relation $R$ is an equivalence relation if and only if the comparison between two finite limits built from $R$ is an isomorphism, see the proof of Proposition 8 in \cite{B-3}.
\endproof
 Let us begin by the following result which is a simple adaptation of Proposition 2.4 in \cite{BG2} from the pointed case to the quasi-pointed one:
\begin{prop}\label{key2}
Suppose $\mathbb E$ is quasi-pointed. A monomorphism $v: U\into T$ with $U\in Pt_0\mathbb E$ has a normalizer in the sense of Definition \ref{orm} if and only if the monomorphism $En v:U \into En T=KJ(T)$ admits a $K$-pre-cartesian monomorphism above it.
This monomorphism $v: U\into T$ is isomorphic to its normalizer if and only if the associated pre-cartesian monomorphism is $\P$-invertible.
\end{prop}
\proof
The full proof is given in \cite{BSM} and mimicks exactly the proof of Proposition 2.4 for pointed categories in \cite{BG2}. It is why, here, we shall only describe the two induced constructions.

Suppose $v$ has a normalizer $(u,R_v)$,  the following right hand side map $(w,\tilde w, v)$ in $Kt\EE$  with $\tilde w=(w.d_0^R,w.d_1^R)$ is pre-cartesian above $En v$:
$$\xymatrix@=20pt{
 {U\;} \ar[d]_{(0,1_U)} \ar@(u,u)@{>->}[rrr]^{En v}\ar@{=}[rr] && {U\;} \ar[d]_{(0,u)} \ar@{>->}[r]^{En v} & En T \ar[d]^{\bar{\epsilon}_T}\\
 {U\times U\;}\ar[d]_{p_0^U} \ar@(u,u)[rrr]_<<<<<<<<<<{v\times v} \ar@{>->}[rr]^{\bar u} && {R_v\;} \ar[d]_{d_0^R} \ar@{>->}[r]^{\tilde w}  & T\times T \ar[d]_{p_0^T}\\
{U\;}\ar@<-1ex>[u]_{s_0^U} \ar@(d,d)[rrr]_{v} \ar@{>->}[rr]_{u} && {X\;} \ar@{>->}[r]_{w} \ar@<-1ex>[u]_{s_0^R} & T \ar@<-1ex>[u]_{s_0^T}
  }
$$
When $v$ is isomorphic to its normalizer, the map $w$ is an isomorphism and this pre-cartesian map is $\P$-invertible.

Conversely, suppose there is a pre-cartesian morphism $(w,\tilde w,Env)$ above $Env$:
$$\xymatrix@=20pt{
 {U\;} \ar[d]_{(0,k_{d_0})} \ar@(u,u)@{>->}[rrr]^{En v} && {U\;} \ar[d]_{k_{d_0}} \ar@{>->}[r]^{En v} & En T \ar[d]^{\bar{\epsilon}_T}\\
 {R[d_0]\;}\ar[d]_{p_0} \ar@(u,u)[rrr]^<<<<<<<<<<<<{(\delta_1.p_0,\delta_1.p_1)} \ar@{.>}[rr]^{d_2} && {R_v\;} \ar[d]_{d_0} \ar@{>->}[r]^{\tilde w}  & T\times T \ar[d]_{p_0^T}\\
{R_v\;}\ar@<-1ex>[u]_{s_0} \ar@(d,d)[rrr]_{\delta_1} \ar@{.>}[rr]_{d_1} && {X\;} \ar@{>->}[r]_{w} \ar@<-1ex>[u]_{s_0} & T \ar@<-1ex>[u]_{s_0^T}
  }
$$
Then the map $\tilde w$ is necessarily of the form $(w.d_0,\delta_1)$ with $\delta_1: R_v\to T$ such that $\delta_1.s_0=w$ (from $\tilde w.s_0=s_0^T.w=(w,w)$) and $\delta_1.k_{d_0}=v$ (from $\tilde w.k_{d_0}=\bar{\epsilon}_T.Env$). The universal property of this pre-cartesian map applied to the morphism $(\delta_1,(\delta_1.p_0,\delta_1.p_1), En v)$ induces a factorization $(d_1,d_2)$ which completes the equivalence relation structure on $R_v$ and makes $u=d_1.k_{d_0}: U\into X$ normal to $R_v$ and such that $w.u=\delta_1.k_{d_0}=v$. This makes $(u,R_v)$ the normalizer of $v$. Saying that the pre-cartesian morphism $(w,\tilde w,Env)$ is $\P$-invertible is  saying that $w$ is invertible and that $v$ is isomorphic to its normalizer.
\endproof
The following lemma is technical and straightforward:
\begin{lemma}\label{tech}
Let $K: \mathbb E\to \mathbb F$ be a left exact functor between finitely complete categories which creates pullbacks. If $u :U\into T$ is a  monomorphism in $\mathbb E$, and $t:T'\to T$ any morphism such that $K(t)$ is a monomorphism and there is map $w$ satisfying $k(u)=K(t).w$, then there is a pullback $\bar u: \bar U\into T'$ of $u$ along $t$ such that $K(\bar u)=w$. If moreover $u$ is $K$-pre-cartesian, so is $\bar u$. In particular, if $u :U\into T$ is a monomorphic $K$-pre-cartesian map in $\mathbb E$ and  $g: W \to T$ a map such that $K(g)$ is an identity map, then there exits a pullback $\bar u$ of $u$ along $g$ such that $K(\bar u)=K(u)$ and $\bar u$ is $K$-pre-cartesian.
\end{lemma}
\begin{prop}\label{key1}
Suppose $\mathbb E$ is quasi-pointed. The following conditions are equivalent:\\
1) the functor $K: Kt\mathbb E\to E$ is pre-fibrant on monomorphisms\\
2) any monomorphism $u: U\into T$ with $U\in Pt_0\mathbb E$ has a normalizer.
\end{prop}
\proof
We have 1) $\Rightarrow$ 2) by the previous proposition. As for the converse, first notice that any split exact sequence as on the left hand side vertical diagram can be embedded in some $J(T)$ in the following way:
$$\xymatrix@=20pt{
 {K_f;} \ar[d]_{k_f}\ar@{>->}[rr]^{En k_f} && {En X\;} \ar[d]_{(0,\epsilon_X)}  \ar@{=}[r] & En X \ar[d]^{\bar{\epsilon}_X}\\
 {X\;}\ar[d]_{f}  \ar@{>->}[rr]^{(f,1)} && {Y\times X\;} \ar[d]_{p_0^Y} \ar@{>->}[r]^{s\times X}  & X\times X \ar[d]_{p_0^X}\\
{Y\;}\ar@<-1ex>[u]_{s} \ar@{=}[rr]  && {Y\;}\ar@{>->}[r]_{s}  \ar@<-1ex>[u]_{(1,s)} & X \ar@<-1ex>[u]_{s_0^X}
  }
$$
Given any monomorphism $m: U\into K_f$, in presence of 2), the monomorphism $En k_f.m: U\into En X=K.J(X)$ has a normalizer which means that there is $K$-pre-cartesian map $(\beta,\alpha,En k_f.m)$ above $En k_f.m$ again by the previous proposition. According to the previous lemma there is a pullback $(\bar{\beta},\bar{\alpha},m)$ of this map along $(s,(s.f,1),En k_f)$ which is necessarily $K$-pre-cartesian above $m$.
\endproof 
In the pointed case, we recover the Theorem 2.8 of \cite{BG2}. The category $Set$ of sets is a non-pointed example satisfying our assumption since the functor $K$, being a terminal functor, is trivially pre-fibrant on monomorphisms. The normalizer of any initial map $\alpha_T: \emptyset\into T$ is nothing but $(\alpha_T,\nabla_T)$, so that any  $\alpha_T$ is isomorphic to its normalizer.

\subsection{Quasi-pointed protomodular categories}\label{mmm}

The protomodular context provides us with a much sharper observation:
\begin{theo}\label{lizer}
	Suppose $\mathbb E$ is quasi-pointed. When $\mathbb E$ is $\P$-decomposable, the functor $K$ is pre-fibrant on monomorphisms.\\ 
	When moreover $\mathbb E$ is protomodular, the following conditions are equivalent:\\
	1) the category $\mathbb E$ is $\P$-decomposable\\
	2) the functor $K$ is pre-fibrant on monomorphisms\\
	3) in the category $\mathbb E$ any monomorphism $u: U\into T$  has a normalizer.\\
	The category $\mathbb E$ is then necessarily stably $\P$-decomposable.
\end{theo}
\proof
Suppose 1). Let $(a,b):A\rightleftarrows B$ be a split epimorphism and $v: U\rightarrowtail K[a]$ any monomorphism. Let us consider the extremal decomposition of the monomorphism $(\alpha_B,v.k_a)$ where $k_a: K[a]\rightarrowtail A$ is the kernel  of $a$:
$$
\xymatrix@=20pt{
	{U\;} \ar[d] \ar@(u,u)@{>->}[rr]^{v.k_a} \ar@{>->}[r]^{\bar u} & {\bar A\;} \ar[d]_{\bar a} \ar@{>->}[r]^{\alpha}  & A \ar[d]_{a}\\
	{0\;} \ar@<-1ex>[u] \ar@(d,d)@{>->}[rr]_{\alpha_B} \ar@{>->}[r]_{\alpha_{\bar B}} & {\bar B\;} \ar@{>->}[r]_{\beta} \ar@<-1ex>[u]_{\bar b} & B \ar@<-1ex>[u]_{b}
}
$$
The left hand side square being a pullback, the map $\bar u$ is a kernel of $\bar a$. Let us show that the following left hand side monomorphism $(\beta,\alpha,v)$ in $Kt\mathbb E$ is $K$-pre-cartesian above $v$:
$$
\xymatrix@=20pt{
	{U\;} \ar@{=}[d] \ar@{=}[r] & {U\;} \ar[d]_{k_f} \ar@(u,u)@{>->}[rr]^{v} & {U\;} \ar[d]_{\bar u} \ar@{>->}[r]^{v} & K[a] \ar[d]^{k_a}\\
	{U\;} \ar[d] \ar@{>->}[r]_{k_f} & {X\;}\ar[d]_{f} \ar@(u,u)[rr]_>>>>>>{x} \ar@{.>}[r]^{\bar x} & {\bar A\;} \ar[d]_{\bar a} \ar@{>->}[r]^{\alpha}  & A \ar[d]_{a}\\
	{0\;}\ar@<-1ex>[u] \ar@{>->}[r]_{\alpha_Y} & {Y\;}\ar@<-1ex>[u]_{s} \ar@(d,d)[rr]_{y} \ar@{.>}[r]_{\bar y} & {\bar B\;} \ar@{>->}[r]_{\beta} \ar@<-1ex>[u]_{\bar b} & B \ar@<-1ex>[u]_{b}
}
$$
So consider any map $(y,x,v)$ in $Kt\mathbb E$. Now complete the diagram in $Kt\mathbb E$ by the map $(\alpha_Y,k_f,1_U)$ on the left hand side; it produces a decomposition of the map $(\alpha_B,v.k_a)$ in $Pt\mathbb E$, whence the dotted factorization $(\bar y,\bar x)$ such that (among other things) $\bar x.k_f=\bar u$, which shows that the factorization of $(\bar y,\bar x)$ at the level of the kernels is $1_U$. Accordingly the map $(\bar y,\bar x,1_U)$ is the required factorization in $Kt\mathbb E$. Whence 2. 

Suppose 2). Then $\mathbb E$ has normalizer for any monomorphism $u: U\into X$ with $U\in Pt_0\mathbb E$ by the previous proposition. When $\mathbb E$ is protomodular, then it has normalizer for any monomorphism by Lemma \ref{rose}. Whence 3). And it is clear that 3) implies 2), again by the previous proposition.

Let us check 2) $\Rightarrow$ 1). First it is easy to check, by Lemma \ref{tech} that the $K$-pre-cartesian maps above the monomorphisms are necessarily monomorphic. Consider any monomorphism $(y,x):(f',s')\rightarrow (f,s)$ in $Pt\mathbb E$. Complete the diagram by the kernels and the factorization $K(x)$, then take the $K$-pre-cartesian map $(\bar y,\bar x,K(x))$ above this monomorphism $K(x)$:
$$
\xymatrix@=20pt{
	{K[f']\;} \ar[d]_{k_{f'}} \ar@{=}[r] \ar@(u,u)@{>->}[rr]^{K(x)} & {K[f']\;} \ar[d]_{k_{\bar f}}\ar@{>->}[r]^{K(x)} & {K[f]\;} \ar[d]^{k_{f}}\\
	{X'\;} \ar[d]_{f'} \ar@(u,u)@{>->}[rr]_>>>>>>{x} \ar@{>.>}[r]^{\underline x} & {\bar X\;} \ar[d]_{\bar f} \ar@{>->}[r]^{\bar x}  & X \ar[d]_{f}\\
	{Y'\;} \ar@<-1ex>[u]_{s'} \ar@(d,d)@{>->}[rr]_{y} \ar@{>.>}[r]_{\underline y} & {\bar Y\;} \ar@{>->}[r]_{\bar y} \ar@<-1ex>[u]_{\bar s} & Y \ar@<-1ex>[u]_{s}
}
$$
It determines a factorization $(\underline y,\underline x,1_{K[f']})$. Since $\mathbb E$ is protomodular, the isomorphic factorization $1_{K[f']}$ at the level of kernels implies that the map $(\underline y,\underline x)$ in $Pt\mathbb E$ is underlying a pullback, namely that this map is $\P_{\mathbb E}$-cartesian. 

It remains to show that the decomposition $(y,x)=(\bar y,\bar x).(\underline y,\underline x)$ is extremal. So, consider another decomposition of $(y,x)=(\bar b,\bar a).(\underline b,\underline a)$ with $(\underline b,\underline a)$ monomorphic and $\P_{\mathbb E}$-cartesian; this implies that the map $\underline a.k_{f'}$ is a kernel of $g$. Complete the following diagram with the kernels:
$$
\xymatrix@=20pt{
	{K[f']\;} \ar[d]_{k_{f'}} \ar@{=}[r] & {K[f']\;} \ar@(u,u)@{>->}[rr]^{K(x)}  \ar[d]_{\underline a.k_{f}} \ar@{=}[r] & {K[f']\;} \ar[d]_{k_{\bar f}}\ar@{>->}[r]^{K(x)} & {K[f]\;} \ar[d]^{k_{f}}\\
	{X'\;} \ar[d]_{f'}  \ar@{>->}[r]^{\underline a} & {A\;} \ar[d]_{g} \ar@(u,u)[rr]_>>>>>>{\bar a} \ar@{.>}[r]^{\alpha} & {\bar X\;} \ar[d]_{\bar f} \ar@{>->}[r]^{\bar x}  & X \ar[d]_{f}\\
	{Y'\;} \ar@<-1ex>[u]_{s'}  \ar@{>->}[r]_{\underline b} & {B\;} \ar@<-1ex>[u]_{t} \ar@(d,d)[rr]_{\bar b} \ar@{.>}[r]_{\beta} & {\bar Y\;} \ar@{>->}[r]_{\bar y} \ar@<-1ex>[u]_{\bar s} & Y \ar@<-1ex>[u]_{s}
}
$$
The map $(\bar b,\bar a,K(x))$ in $Kt\mathbb E$ gives a unique factorization $(\beta,\alpha,1_{K[f']})$ through the $K$-pre-cartesian map $(\bar y,\bar x,K(x))$. The map $(\beta,\alpha)$ is actually the desired factorization in $Pt\mathbb E$, so that $\mathbb E$ is $\P$-decomposable. The fact that $\mathbb E$ is stably $\P$-decomposable is a consequence of the fact that the $K$-pre-cartesian morphisms are necessarily stable under pullbacks along $\P_{\mathbb E}$-cartesian morphisms, again according to Lemma \ref{tech}, since any $\P_{\mathbb E}$-cartesian morphism can be extented into a morphism in $Kt\mathbb E$ whose image by $K$ is an identity map.
\endproof

\begin{coro}\label{truco}
	Let $\mathbb C$ be quasi-pointed, protomodular and $\P$-decomposable. Then:\\
	1) a morphism $(y,x):(\bar f,\bar s)\rightarrow (f,s)$ is $\P$-cartesian if and only if it is $K$-invertible;\\
	2) a monomorphism $(y,x):(f',s')\rightarrowtail (f,s)$ is a $\P$-outsider if and only if it determines a $K$-pre-cartesian monomorphism above $K(x)$ in $Kt\mathbb E$:
	$$
	\xymatrix@=20pt{
		{K[\bar f]\;} \ar[d]_{k_{\bar f}}\ar@{>->}[r]^{K(x)} & {K[f]\;} \ar[d]^{k_{f}}\\
		{\bar X\;} \ar[d]_{\bar f} \ar@{>->}[r]^{x}  & X \ar[d]_{f}\\
		{\bar Y\;} \ar@{>->}[r]_{y} \ar@<-1ex>[u]_{\bar s} & Y \ar@<-1ex>[u]_{s}
	}
	$$
	3) a monomorphism $v: U\into T$ is normal in $\CC$ iff and only if  the $K$-pre-cartesian morphism associated with $En v: U\into En T$ is $\P$-invertible.
\end{coro}
\proof
The first point is a classical characterization of pullbacks in $Pt\CC$. According to the proof of [2) $\rightarrow$ 1)] in the previous theorem, the monomorphism $(y,x)\in Pt\CC$ is a $P$-outsider if and only if the map $\underline y$ is $1_Y$, namely if and only if $(y,x)$ is a $\P$-outsider; whence 2). In our protomodular context, a monomorphism  is normal if and only if it is isomorphic to its normalizer; then 3) is a consequence of Proposition \ref{key2}.
\endproof
\begin{coro}\label{xgr}
	Let $(\;)_0:Grd\rightarrow Set$ be the forgetful functor from groupoids to sets. Then any fibre $Grd_Y$ is stably $\P$-decomposable.
\end{coro}
\proof
We recalled that any fibre $Grd_Y$ is quasi-pointed and protomodular. A subobject $\underline u_1: \underline U_1\hookrightarrow \underline X_1$ is normal in $Grd_Y$ if and only if, for any arrow  $\phi:y\to y'$ in $\underline X_1$ and any endomap $\tau$ on $y$ in $\underline U_1$, the endomap $\phi.\tau.\phi^{-1}$ is in $\underline U_1$ \cite{B15}. The normalizer of any subobject $\underline v_1:\underline U_1\into \underline T_1$ in $Grd_Y$ is then defined by the subset $\underline X_1$ of those arrows $\phi:y\to y'$ of $\underline T_1$ which are such that, for any endomap $\tau$ on $y$ in $\underline U_1$, the endomap $\phi.\tau.\phi^{-1}$ is in $\underline U_1$ and, for any endomap $\theta$ on $y'$ in $\underline U_1$, the endomap $\phi^{-1}.\theta.\phi$ is in $\underline U_1$.
\endproof

In Section \ref{groupe}, the construction of the $\P$-decomposition in the category $Gp$ of groups needed the universal quantifier. It is  why it is necessary to add the assumption of cartesian closedness to show that the category $Gp\EE$ of internal groups in $\EE$ has normalizers, see \cite{BG2}. Of course, the  universal quantifier is also used in the previous corollary. It is showed in the preprint \cite{BSM} that, in the same way as for  $Gp\EE$, any fibre $Grd_Y\mathbb E$ is stably $\P$-decomposable,  provided that the ground category $\EE$ is locally cartesian closed (which is true when  $\EE$ is a topos, see \cite{Jo}).

\subsection{The fibers $Cat_Y\EE$}

Let $(\;)_0: Cat \to Set$ be the fibration associating to any category its set of objects. The category $Mon$ of monoids is the fiber above $1$. Any other fiber $Cat_Y$ is quasi-pointed. Call Schreier split epimorphism any split epimorphism $(\underline f_1,\underline s_1)$ in this fiber such that $\underline f_1$ is a split fibration or, in other words, any split epimorphism $(\underline f_1,\underline s_1)$, where the splitting $\underline s_1$ is a fibrant splitting. This class $\Sigma_Y$ of split epimorphisms behaves in $Cat_Y$ exactly as the class of Schreier split epimorphisms in $Mon$ and determines a subfibration $\P_{\Sigma_Y}$ of $\P_{Cat_Y}$, see \cite{85}. On the model of $Mon$, see \cite{BSM}, we get: 
\begin{prop}
	Any fiber $Cat_Y$ is stably $\P_{\Sigma_Y}$-decomposable. When $\EE$ is locally cartesian closed, it is the case for any fiber $Cat_Y\EE$.
\end{prop}
\proof
Given any subobject in $Pt(Cat_Y)$ between split epimorphims with fibrant splittings:
$$
\xymatrix@=20pt{
	{\mathbb X'\;}\ar[d]_{F'} \ar@<2ex>@{ >->}[rr]_>>>>>>{} \ar@{>.>}[r]^{} & {\bar{\mathbb X}\;} \ar[d]_{\bar F} \ar@{>.>}[r]^{}  & \mathbb X \ar[d]_{F}\\
	{\mathbb Y'\;}\ar@<-1ex>[u]_{S'} \ar@<-2ex>@{ >->}[rr]_{} \ar@{>.>}[r]_{} & {\bar{\mathbb Y}\;} \ar@{>.>}[r]_{} \ar@<-1ex>[u]_{\bar S} & \mathbb Y \ar@<-1ex>[u]_{S}
}
$$
the morphisms of the subcategory $\bar{\mathbb Y}$ are those morphisms  $\psi\in \mathbb Y$ such that $q(S(\psi).u)\in \mathbb X', \;\forall u\in Ker F'$ while the morphisms of the subcategory $\bar{\mathbb X}$ those morphisms $\phi\in \mathbb X$ such that $F(\phi)\in \bar{\mathbb Y}, \;{\rm and}\; q(\phi)\in \mathbb X'$. From this construction, checking the universal property and the stablitily under pullback along $\P_{\Sigma}$-cartesian morphisms is straightforward.
In the internal context $Cat_Y\EE$, the proof is detailed in \cite{BSM}; and again the assumption of the locally cartesian closedness is needed by the presence of the universal quantifier in the definition of the subcategory $\bar{\mathbb Y}$.  
\endproof

\section{$\P$-decomposable category and  additive setting}\label{adset}

Let us recall that a pointed category $\mathbb A$ is additive if and only if any object $X$ is endowed with a natural internal group structure which is then necessarily commutative and that a first kind of "non-pointed additive" setting was introduced with the following definition and proposition:
\begin{defi}\cite{Jo2}
	A category $\mathbb A$ is said to be a naturally Mal'tsev category when any object $X$ is endowed with a natural Mal'tsev operation (namely ternary operation $p_X:X\times X\times X \to X$ satisfying $p_X(x,y,y)=x=p_X(y,y,x)$)  which is then necessarily associative and commutative.
\end{defi}
\begin{prop}\cite{B-3}
	A category $\mathbb A$ is a naturally Mal'tsev one if and only if any fibre $Pt_Y\mathbb A$ is additive.
\end{prop}
In Proposition \ref{thinvert}, we observed that in a $\P$-decomposable  category $\EE$ any $\P$-invertible monomorphism is a $\P$-outsider. Here we shall investigate the inverse implication and show  that it is strongly related with this non-pointed additive setting.

\subsection{$\P$-invertible vs $\P$-outsider}  

\begin{prop}\label{outinv}
	Consider the following conditions:\\
	1) the category $\mathbb E$ is $\P$-decomposable and the $\P$-outsiders coincide with the $\P$-invertibles;\\
	2) the category $\mathbb E$ is a Mal'tsev one and any base-change $y^*$ with respect to the fibration $\P_{\mathbb E}$ along a monomorphism $y$ is pre-fibrant on subobjects;\\
	3) the category $\mathbb E$ is a Mal'tsev one and, for any monomorphism $m: Y'\into Y$ and any equivalence relation $R$ on $Y'$, there is a largest equivalence relation $\Sigma$ on $Y$ among those ones which are such that $m^{-1}(S)$ is $R$ and the induced monomorphism $R\into \Sigma$ is fibrant;\\
	4) the category $\mathbb E$ is a Mal'tsev one and, for any monomorphism $m: Y'\into Y$, there is a largest equivalence relation $R_m$ to which $m$ is normal;\\
	$\alpha$) $\mathbb E$ is a Mal'tsev category with normalizers and, for any monomorphism $m: Y'\into Y$, there is a largest equivalence relation $R_m$ to which $m$ is normal.\\
	Then we get 1)$\Rightarrow$ 2)$\Rightarrow$ 3) $\Rightarrow$ 4) and 1)$\Rightarrow$ $\alpha$)$\Rightarrow$ 4).\\
	If, moreover, $\mathbb E$ is protomodular we get: 1)$\iff$ 2)  and $\alpha$)$\iff$ 4) $\Rightarrow$ 5), with:\\  5) $\EE$ is a naturally Mal'tsev category.
\end{prop}
\proof
Suppose 1). By Corollary \ref{malc2}, we know that $\mathbb E$ is a Mal'tsev category. Given any left hand side diagram where the square is a pullback:
$$ \xymatrix@=8pt{
	{X''\;} \ar@{>->}[rrr]^{\xi} \ar@<-1ex>[dd]_{f''} \ar@{>->}@(u,u)[rrrrrr]^{x} &&& {X'\;} \ar@<-1ex>[dd]_{f'} \ar@{>->}[rrr]^{x'} &&& X \ar@<-1ex>[dd]_{f}&&   {X''\;} \ar@{>->}[rrr]^{x''} \ar@<-1ex>[dd]_{f''} \ar@{>->}@(u,u)[rrrrrr]^{x} &&& {\bar X\;} \ar@<-1ex>[dd]_{\bar f} \ar@{>->}[rrr]^{\bar{\xi}} &&& X \ar@<-1ex>[dd]_{f}\\
	&  &&&  \\
	{Y'\;} \ar@{=}[rrr]  \ar[uu]_{s''} &&& Y'  \ar[uu]_{s'} \ar@{>->}[rrr]_{y} &&& Y \ar[uu]_{s}&&   {Y'\;} \ar@{>->}[rrr]_{y}  \ar[uu]_{s''} &&& Y  \ar[uu]_{\bar s} \ar@{=}[rrr] &&& Y \ar[uu]_{s}
}
$$
the extremal decomposition of $(y,x'')$ on the right hand side produces the desired $m^*$-pre-cartesian map above $\xi$. Whence 2).

Suppose 2). Then consider the $y^*$-pre-cartesian map associated with the monomorphism $(d_0^R,y.d_1^R)$ in $Pt_{Y'}\EE$ in the left hand side diagram which produces a reflexive relation $\Sigma$ on $Y$ as on the right hand side and a monomorphism $(y,\tilde y): (d_0^R,s_0^R)\into (d_0^{\Sigma},s_0^{\Sigma})$ of split epimorphisms: 
$$ \xymatrix@=8pt{
	{R\;} \ar@{>->}[rrr]^{(d_0^R,y.d_1^R)} \ar@<-1ex>[dd]_{d_0^R} \ar@{>->}@(u,u)[rrrrrr]^{(y.d_0^R,y.d_1^R)} &&& {Y'\times Y\;} \ar@<-1ex>[dd]_{p_0^{Y'}} \ar@{>->}[rrr]^{y\times Y} &&& Y\times Y \ar@<-1ex>[dd]_{p_0^Y} &&   {R\;} \ar@{>->}[rrr]^{\tilde y} \ar@<-1ex>[dd]_{d_0^R} \ar@{>->}@(u,u)[rrrrrr]^{(y.d_0^R,y.d_1^R)} &&& {\Sigma\;} \ar@<-1ex>[dd]_{d_0^{\Sigma}} \ar@{>->}[rrr]^{(d_0^{\Sigma},d_1^{\Sigma})} &&& Y\times Y \ar@<-1ex>[dd]_{p_0^Y}\\
	&  &&&  \\
	{Y'\;} \ar@{=}[rrr]  \ar[uu]_{s_0^R} &&& Y'  \ar[uu]_{(1,y)} \ar@{>->}[rrr]_{y} &&& Y \ar[uu]_{s_0^Y}&&   {Y'\;} \ar@{>->}[rrr]_{y}  \ar[uu]_{s_0^R} &&& Y  \ar[uu]_{s_0^{\Sigma}} \ar@{=}[rrr] &&& Y \ar[uu]_{s_0^Y}
}
$$
This reflexive relation $\Sigma$ is actually an equivalence relation since $\mathbb E$ is a Mal'tsev category and the monomorphism $(y,\tilde y)$ becomes a monomorphism $R\into \Sigma$ of equivalence relations which is fibrant since the square indexed by $0$ is a pullback. Accordingly we get $y^{-1}(\Sigma)=R$. The universal property
of the $y^*$-pre-cartesian map shows that $\Sigma$ is the largest equivalence relation on $Y$ among those $S$ which are such that $y^{-1}(S)$ is $R$ and the induced monomorphism $R\into \Sigma$ is fibrant. Whence 3). The implication 3) $\Rightarrow$ 4) is trivial. The implication 1) $\Rightarrow$ $\alpha$) is a consequence of Proposition \ref{lizer1} and of ($1\Rightarrow 4)$). Again the implication $\alpha$) $\Rightarrow$ 4) is trivial.

Suppose now that $\mathbb E$ is protomodular and 2) is satified. Take any monomorphism in $Pt\EE$ as in the left hand side diagram and take the pullback of  $(f,s)$ as in the middle diagram:
$$
\xymatrix@=20pt{
	{X'\;} \ar@{>->}[r]^{x} \ar@<-1ex>[d]_{f'}  & {X\;} \ar@<-1ex>[d]_{f} & {X'\;} \ar[d]_{f'} \ar@(u,u)@{>->}[rr]^{x} \ar@{>->}[r]^{\bar u} & {\bar X\;} \ar[d]_{\bar f} \ar@{>->}[r]^{\bar w}  & X \ar[d]_{f} & {X'\;} \ar[d]_{f'} \ar@(u,u)@{>->}[rr]^{x} \ar@{>->}[r]^{\check w} & {\bar X'\;} \ar[d]_{\bar f'} \ar@{>->}[r]^{\check  u}  & X \ar[d]_{f}\\
	{Y'\;} \ar@{>->}[r]_{y}  \ar[u]_{s'} & Y  \ar[u]_{s}  & {Y\;} \ar@<-1ex>[u]_{s'}  \ar@{=}[r] & {Y\;} \ar@{>->}[r]_{y} \ar@<-1ex>[u]_{\bar s} & Y \ar@<-1ex>[u]_{s} & {Y\;} \ar@<-1ex>[u]_{s'}  \ar@{>->}[r]_y & {\;} \ar@{=}[r] \ar@<-1ex>[u]_{\bar s} & Y \ar@<-1ex>[u]_{s}
}
$$
This produces a monomorphism in $Pt_Y\EE$ which by 2) produces the right hand side above diagram where its left hand part is a pullback. By Proposition \ref{malc1}, it is an extremal $\P$-decomposition. That any extremal $\P$-decomposition is of this kind implies  that any $\P$-outsider is $\P$-invertible.

Suppose $\EE$ protomodular and $\alpha$). By Proposition \ref{surnorm}, it is equivalent to saying that \emph{$\EE$ is a protomodular category and any monomorphism  is normal}. Then 4) becomes trivial.

Let us show 5). When $\EE$ is protomodular, so is any (pointed) fiber $Pt_Y\EE$. A category is additive if and only if it is pointed protomodular and such that any monomorphism is normal \cite{BB}. 

Suppose $\EE$ is protomodular and such that any monomorphism is normal. Let us show that any monomorphism $m$ in $Pt_Y\EE$ is normal in this fiber. Let $R$ denote the equivalence relation on $X$ to which the monomorphism $m$ is normal in the ground category $\EE$:
$$
\xymatrix@=20pt{
	{X'\;} \ar@{>->}[r]^{m} \ar@<-1ex>[d]_{f'}  & {X\;} \ar@<-1ex>[d]_{f} &  {X'\times X'\;} \ar@<-1ex>[d]_{p_0^{X'}} \ar@<1ex>[d]^{p_1^{X'}} \ar@{>->}[r]^{\tilde m}  & R \ar@<-1ex>[d]_{d_0^R} \ar@<1ex>[d]^{d_1^R} & {R[f']\;} \ar@{ >->}[d]_{}  \ar@{>->}[r]^{} & {R[f] \cap R\;} \ar@{ >->}[d]_{} \\
	{Y\;} \ar@{=}[r]  \ar[u]_{s'} & Y  \ar[u]_{s}  & {X'\;} \ar@{>->}[r]_{m} \ar[u]_{} & X \ar[u]_{} &  {\nabla_{X'}\;}  \ar@{>->}[r]_{(m,\tilde m)} & {R\;} 
}
$$
The monomorphism $m$ lying in the fiber $Pt_Y\EE$ determines a cartesian morphism $R[f']\into R[f]$ above $m$. Whence: $m^{-1}(R[f]\cap R)=m^{-1}(R[f])\cap m^{-1}(R)=R[f']\cap \nabla_{X'}=R[f']$; and the above right hand side pullback in $Grd\EE$. Accordingly the monomorphism $R[f']\into R[f]\cap R$ is fibrant, since so is $\nabla_{X'}\into R$. By $R[f]\cap R\subset R[f]$ the equivalence relation $R[f]\cap R$ lies in $Pt_Y\EE$. So, the monomorphism $m$ of $Pt_Y\EE$ is normal to the equivalence relation $R[f]\cap R$ in $Pt_Y\EE$, and $Pt_Y\EE$ is additive.
\endproof

\begin{coro}
	Let $\EE$ be a quasi-pointed protomodular category $\EE$. The following conditions are equivalent:\\
	1) $\mathbb E$ is $\P$-decomposable and such that any $\P$-outsider monomorphism is $\P$-invertible;\\
	2) the functor $K: Kt\mathbb E \to Pt_0\mathbb E$ is pre-fibrant on monomorphism, and any pre-cartesian monomorphism is $\P$-invertible:\\
	3) any monomorphism in $\mathbb E$ is a normal monomorphism.\\
	In this case, $\EE$ is a naturally Mal'tsev category. In particular, a pointed protomodular category $\EE$ is additive if and only if it is is $\P$-decomposable and such that any $\P$-outsider monomorphism is $\P$-invertible.
\end{coro} 
\proof
By 1) $\iff$ 2) in Theorem \ref{lizer} and by Corollary \ref{truco}, when $\mathbb E$ is quasi-pointed and protomodular we have immediately 1) $\iff$ 2). We noticed that, in the protomodular context, saying that any monomorphism is normal is saying that any monomorphism is isomorphic to its normalizer; so 1) $\iff$ 3) is a consequence of  Theorem \ref{lizer} and of Proposition \ref{key2}. 
The last point is then a straighforward consequence of Condition 5) in the previous proposition.
\endproof
An example of such a quasi-pointed context is given by any fiber $Grd_Y\EE$ in a Mal'tsev category $\EE$, see Corollary 2.8 in \cite{B15}.

\subsection{"Non-pointed additive" settings}   

 The naturally Mal'tsev setting is the largest step of a decreasing scale of four among the "non-pointed additive" settings:
\begin{defi}\cite{B15}
	A category $\mathbb C$ is:\\
	2) antepenessentially affine when any base-change functor is fully faithful;\\
	3) penessentially affine when, in addition, any base-change functor is fibrant on monomorphisms;\\
	4) essentially affine when any base-change functor is an equivalence of categories.	
\end{defi}  
In \cite{B15}, a strict example of each level is given. These three  types of non-pointed additive categories  are necessarily protomodular since a fully faithful functor is necessarily conservative.  They are naturally Mal'tsev categories; and in a penessentially category, any morphism in normal, again see \cite{B15}.
We can now add, in the exact context \cite{Ba}, some precision about the relationship between the above two first levels:
\begin{prop}
	Consider the following conditions:\\
	1) $\mathbb C$ is penessentially affine\\
	2) $\mathbb C$ is an antepenessentially affine, $\P$-decomposable and the $\P$ outsiders coincide with the $\P$-invertibles;\\
	3) $\mathbb C$ is protomodular,  $\P$-decomposable and the $\P$ outsiders coincide with the $\P$-invertibles.\\
	We get  1) $\Rightarrow$ 2) $\Rightarrow$ 3). When, in addition, $\mathbb{C}$ is exact, these three conditions are  equivalent.
\end{prop}
\proof
By definition of a penessentially affine category and by Proposition \ref{outinv}, we have [1)$\Rightarrow$ 2)]. We already noticed that any antepenessentially affine is protomodular, whence [2)$\Rightarrow$ 3)]. Now suppose $\CC$ exact and 3). Any protomodular category is a Mal'tsev one. In a regular (and a fortiori exact) Mal'tsev category, any base-change functor $f^*$ along a regular epimorphism $f$ is fully faithful and fibrant on monomorphisms by Theorem 52 in \cite{BGJ}. It remains to show that, in the exact context, this is the case for any base-change $m^*$ along a monomorphims $m$ as well.

The base-change $m^*$ is conservative since $\mathbb C$ is protomodular, and pre-fibrant on monomorphims by Proposition \ref{outinv}; so, according to Proposition \ref{key6}, it is fibrant on monomorphims. It remains to show that $m^*$ is  fully faithful. This will be the consequence of the following lemma.
\endproof

\begin{lemma}
	Let $U: \EE \to \FF$ be any left exact functor which is fibrant on monomorphism. Then $U$ is "fully faithful on monomorphism". When, in addition, $\EE$ is an exact category and $\FF$ a regular one, then $U$ is fully faithful; it reflects and preserves  the regular epimorphims.
\end{lemma}
\proof
We know by Proposition \ref{key6} that $U$ is conservative. Any left exact conservative functor is faithful. First observe that two objects $X$ and $X'$ with same image $Y$ by $U$ are isomorphic above $1_Y$ in unique way. We have $U(X\times X')=Y\times Y$. Take the fibrant monomorphism $i: W\into X\times X'$ above the diagonal $s_0^Y: Y\into Y\times Y$. Then $U(p_X.i)=1_Y$, so $p_X.i: W\to X$ is an isomorphism, and the conclusion. From that, $U$ is "fully faithful on isomorphisms". Suppose now you have a monomorphism $n: U(X)\into U(X')$. Let $m: W\into X'$ the cartesian  monomorphism above it. Then $U(W)=U(X) $ and by the previous isomorphism, you get a monomorphism above $n$.

Suppose, in addition, $\EE$ exact and $\FF$ regular. Let $\phi: U(X)\onto U(X')$ be any regular epimorphism. Take $S$ the equivalence relation on $X$ above the kernel equivalence relation $R[\phi]$ on $U(X)$ given by Lemma \ref{key6}. This equivalence relation $S$ a quotient $q:X\onto W$ since $\EE$ is exact. Its image $u(q)$ has $R[f]$ as kernel equivalence relation. Since $\EE$ is regular, you get a monomorphic factorizaton $n: U(X')\into U(Y)$ such that $n.\phi=U(q)\; (*)$. Let $m$ be the cartesian monomorphism above it. Denote $\bar m: W\into X$ its pullback along $q$, and $\bar q: W\onto X'$ the induced regular epimorphism. Then $(*)$ makes $U(\bar m)$ invertible; and since $U$ is conservative, $m$ itself is invertible. So, $U(\bar q.\bar m^{-1})=\phi$. This makes $U$ fully faithful. Its reflects the regular epimorphism since $m.f=q$ and  $q$ being a regular epimorphism, the monomorphism $m$ is an isomorphism; so, $f$ is itself a regular epimorphism. 

Finally let us show that $U$ preserves the regular epimorphisms. Let $f:X\onto X'$ be a regular epimorphism in $\EE$. The functor  $U$ preserves $R[f]$, so $U(R[f])$ is an effective equivalence relation in the regular category $\FF$. Denote $\bar q: U(X)\onto Y$ its quotient and $n:Y\into U(X')$ the monomorphic factorization such that $n.\bar q=U(f) \;(*)$. Let $m: W\into X'$ be the cartesian monomorphism above it.
Then the factorization $(*)$ induces a factorization $m.q=f$. Since $f$ is a regular epimorphism, $m$ is an isomorphism. So is $n=U(m)$, and $U(f)$ is a regular epimorphism. 
\endproof
\begin{coro}
	Let $\EE$ be an exact pointed protomodular category. The two following conditions are equivalent:\\
	1) $\EE$ is penessentially affine;\\
	2) any monomorphism is normal.
\end{coro}

\vspace{5mm}
\noindent Univ. Littoral C\^ote d'Opale, UR 2597, LMPA,\\
Laboratoire de Math\'ematiques Pures et Appliqu\'ees Joseph Liouville,\\
F-62100 Calais, France. bourn@univ-littoral.fr

\end{document}